\documentclass[11pt]{amsart}
\usepackage{enumerate}
\usepackage{amssymb,amsthm,amscd}

\newcommand{\op}[1]{\operatorname{#1}}
\newcommand{\volume}{\op{vol}}
\newcommand{\Cal}[1]{\mathcal{#1}}
\renewcommand{\Bbb}[1]{\mathbb{#1}}
\allowdisplaybreaks

\newcommand{\svolball}[2]{{\volume(\underline B_{#2}^{#1})}}
\newcommand{\svolann}[2]{{\volume(\underline A_{#2}^{#1})}}
\newcommand{\sball}[2]{{\underline B_{#2}^{#1}}}
\newcommand{\svolsp}[2]{{\volume(\partial \underline B_{#2}^{#1})}}
\newtheorem*{thmA}{Theorem A}
\newtheorem*{thmB}{Theorem B}
\newtheorem*{thmC}{Theorem C}
\newtheorem*{thmD}{Theorem D}
\newtheorem*{thmE}{Theorem E}

\newtheorem{theorem}{Theorem}[section]
\newtheorem{corollary}[theorem]{Corollary}
\newtheorem{lemma}[theorem]{Lemma}

\theoremstyle{remark}
\newtheorem{remark}[theorem]{Remark}
\newtheorem{conjecture}[theorem]{Conjecture}

\begin{document}

\title[Quantitative Volume Rigidity]{Quantitative Volume Space Form Rigidity Under Lower Ricci Curvature Bound I}

\author{Lina Chen}
\address{Department of mathematics, Nanjing University, Nanjing China}
\email{chenlina\_mail@163.com}
\thanks{The first author is supported partially by a research fund from Capital Normal University.}

\author{Xiaochun Rong}
\address{Mathematics Department, Rutgers University
	New Brunswick, NJ 08903 USA}

\email{rong@math.rutgers.edu}
\thanks{The second author is supported partially by NSF
	Grant DMS 0203164 and by a research fund from Capital Normal University.}
	
\author{Shicheng Xu}
\address{Mathematics Department and Academy for Multidisciplinary Studies, Capital Normal University, Beijing China}
\email{shichengxu@gmail.com}
\thanks{The third author is supported partially by NSFC Grant 11821101, 11871349 and by Youth Innovative Research Team of Capital Normal University.}
\date{}

\begin{abstract} Let $M$ be a compact $n$-manifold of $\op{Ric}_M\ge (n-1)H$ ($H$ is a constant).
	We are concerned with the following space form rigidity: $M$ is isometric to a space
	form of constant curvature $H$ under either of the following conditions:
	
	(i) There is $\rho>0$ such that for any $x\in M$, the open $\rho$-ball at $x^*$ in the
	(local) Riemannian universal covering space, $(U^*_\rho,x^*)\to (B_\rho(x),x)$,
	has the maximal volume i.e., the volume of a $\rho$-ball in the simply connected $n$-space
	form of curvature $H$.
	
	(ii) For $H=-1$, the volume entropy of $M$ is maximal i.e. $n-1$ ([LW1]).
	
	\noindent The main results of this paper are quantitative space form rigidity i.e.,
	statements that $M$ is diffeomorphic and close in the Gromov-Hausdorff topology to
	a space form of constant curvature $H$, if $M$ almost satisfies, under some additional condition,
	the above maximal volume condition. For $H=1$, the quantitative spherical space form rigidity
	improves and generalizes the diffeomorphic sphere theorem in [CC2].
\end{abstract}

\maketitle

\setcounter{section}{-1}
\section{Introduction}
Let $M$ be a compact $n$-manifold of $\op{Ric}_M\ge (n-1)H$, $H$ is a constant. The goal of this
paper is to establish quantitative version for two space form rigidity under lower Ricci curvature bound
(see Theorem 0.1 and 0.3).  Our quantitative version has two components: it includes a rigidity and
reveals a diffeomorphism stability.  This work is based on, among other things,
the work of Cheeger-Colding (\cite{Ch}, \cite{Co1,Co2}, \cite{CC1,CC2}).

The first one is essentially the rigidity part of Bishop volume comparison. For our purpose
(see Conjecture 0.15), we formulate it as follows. For a metric ball $B_r(x)$ on a manifold $M$,
we will call $B_r(x^*)$ a local rewinding of $B_r(x)$ and the volume, $\op{vol}(B_r(x^*))$,
the local rewinding volume of $B_r(x)$, where $\pi^*: (U^*_\rho,x^*)\to (B_\rho(x),x)$ is the
(incomplete) Riemannian universal covering space. Similarly, if $\pi: (\tilde M,\tilde x)\to (M,x)$ is
a Riemannian universal cover, we call $B_r(\tilde x)$ a global rewinding of $B_r(x)$.

\begin{theorem}[Maximal local rewinding volume rigidity]
	Let $M$ be a compact $n$-manifold of $\op{Ric}_M\ge (n-1)H$.
	If there is $\rho>0$ such that for any $x\in M$, the local rewinding volume
	$\op{vol}(B_\rho(x^*))=\svolball{H}{\rho}$, then $M$ is isometric to a space form of curvature $H$,  where
	$\sball{H}{\rho}$ denotes a $\rho$-ball in the simply
	connected $n$-space form of constant curvature $H$.
\end{theorem}

For $H\ge 0$, $M$ in Theorem 0.1 may have an arbitrarily small volume i.e., collapsed.
For $H=1$, Theorem 0.1 includes the maximal volume rigidity: if a complete
$n$-manifold $M$ of $\op{Ric}_M\ge n-1$ achieves the maximal volume (when $\rho=\pi$) i.e.,
the volume of unit sphere, then $M$ is isometric to $S^n_1$.

A quantitative maximal volume rigidity is the following
sphere theorem:

\begin{theorem}[\cite{CC2}]
	There exists a constant $\epsilon(n)>0$
	such that for any $0\le \epsilon<\epsilon(n)$, if a compact $n$-manifold $M$ satisfies
	$$\op{Ric}_M\ge n-1,\quad \frac{\volume(M)}{\volume(S^n_1)}\ge 1-\epsilon,$$
	then $M$ is diffeomorphic to the unit sphere, $S^n_1$, by a $\Psi(\epsilon|n)$-isometry
	(i.e., a diffeomorphism with a distance distortion at most $\Psi(\epsilon|n)$),
	where $\Psi(\epsilon|n)\to 0$ as $\epsilon\to 0$ while $n$ is fixed.
\end{theorem}

A homeomorphism in Theorem 0.2 was first obtained in \cite{Pe1}, a $\Psi(\epsilon|n)$-closeness
was established in \cite{Co1}, and Theorem 0.2 was proved in \cite{CC2} via the Reifenberg's method.
Note that Theorem 0.2 implies the maximal volume rigidity.

The other space form rigidity result is the Ledrappier-Wang's maximal volume entropy
rigidity (see Theorem 0.3). The volume entropy of a compact manifold $M$ is defined by
$$h(M)=\lim_{R\to \infty}
\frac{\ln(\volume(B_R(\tilde p)))}R,\qquad \tilde p\in \tilde M$$
(for the existence of the limit, see \cite{Ma}), where $\tilde M$ denotes the Riemannian
universal covering space of $M$. By Bishop volume comparison, for
any compact $n$-manifold $M$ of $\op{Ric}_M\ge -(n-1)$, $h(M)\le n-1$,
which equals to the volume entropy of any hyperbolic $n$-manifold.

\begin{theorem}[Maximal volume entropy rigidity \cite{LW1}]
If a compact $n$-manifold $M$ of $\op{Ric}_M\ge -(n-1)$
achieves the maximal volume entropy i.e., $h(M)=n-1$,  then $M$ is
isometric to a hyperbolic manifold.
\end{theorem}

We now begin to state our quantitative version for Theorem 0.1 with
respect to local rewinding volume and normalized $H=\pm 1$ and $0$ respectively,
starting with $H=1$.

\begin{thmA}
	Given $n, \rho, v>0$, there exists a constant
	$\epsilon(n,\rho,v)>0$ such that for any $0\le \epsilon<\epsilon(n,\rho,v)$,
	if a compact $n$-manifold $M$ satisfies
	$$\op{Ric}_M\ge n-1,\quad\op{vol}(\tilde M)
	\ge v,\quad \frac{\volume
		(B_\rho(x^*))}{\svolball{1}{\rho}}\ge 1-\epsilon,\quad\forall\, x\in M,$$
	then $M$ is diffeomorphic to a spherical space form by a $\Psi(\epsilon|n,\rho,v)$-isometry,
	where $\op{vol}(B_\rho(x^*))$ denotes the local rewinding
	volume of $B_\rho(x)$.
\end{thmA}

Theorem A generalizes and improves Theorem 0.2, see Remark 0.7.
For $H=-1$, we have

\begin{thmB}
	Given $n, \rho, d, v>0$, there exists
	$\epsilon(n,\rho,d,v)>0$ such that for any $0\le \epsilon<\epsilon(n,\rho,v,d)$,
	if a compact $n$-manifold $M$ ($\tilde p\in \tilde M$) satisfies that for all $x\in M$,
	$$\op{Ric}_M\ge -(n-1),\,\op{diam}(M)\le d,\,\op{vol}(B_1(\tilde p))\ge v,\,\frac{\volume
		(B_\rho(x^*))}{\svolball{-1}{\rho}}\ge 1-\epsilon,$$
	then $M$ is diffeomorphic to a hyperbolic manifold by a $\Psi(\epsilon|n,\rho,d,v)$-isometry.
\end{thmB}

Note that Theorem B does not hold if one removes a bound on diameter; there is a sequence
of compact $n$-manifolds $M_i$ ($n\ge 4$) of negative pinched sectional
curvature $-1\le \op{sec}_{M_i}\le -1+\epsilon_i$ and
$\epsilon_i\to 0$ ($\op{diam}(M_i)\to \infty$), but $M_i$ admits
no hyperbolic metric (\cite{GT}).  On the other hand, given any $\rho, \epsilon>0$,
it is clear that for $i$ large, $\frac{\op{vol}(B_\rho(\tilde x_i))}{\svolball{-1}
	{\rho}}\ge 1-\epsilon$ for any $\tilde x_i\in \tilde M_i$.

For $H=0$, because of the splitting theorem of Cheeger-Gromoll (\cite{CG}) we actually
prove a rigidity result.

\begin{thmC}
	Given $n, \rho, v>0$, there exists $\epsilon=\epsilon(n,\rho,v)>0$
	such that if a compact $n$-manifold $M$ ($\tilde p\in \tilde M$) satisfies
	$$\op{Ric}_M\ge 0,\, \op{diam}(M)\le 1,\, \op{vol}(B_1(\tilde p))\ge v,\,
	\frac{\volume(B_\rho(x^*))}
	{\svolball{0}{\rho}}\ge 1-\epsilon,
	\,\forall\, x\in M,$$
	then $M$ is isometric to a flat manifold.
\end{thmC}

A quantitative version of Theorem C is the following.

\begin{theorem}
	Given $n, \rho, v>0$, there exist $\delta(n,\rho,v),
	\epsilon(n,\rho,v)>0$ such that for any $0<\delta<\delta(n,\rho,v)$,
	if a compact $n$-manifold $M$ satisfies that for all $x\in M$,
	$$\op{Ric}_M\ge -(n-1)\delta,\,1\ge \op{diam}(M),\,\op{vol}
	(M)\ge v,\,\frac{\volume (B_\rho(x^*))}
	{\svolball{0}{\rho}}\ge 1-\epsilon(n,\rho,v),$$
	then $M$ is diffeomorphic to a flat manifold by a $\Psi(\delta|n,\rho,v)$-isometry.
\end{theorem}

Note that unlike Theorem A-C, Theorem 0.4 does not hold if one relaxes the condition, `$\op{vol}(M)\ge v$', to
`$\op{vol}(B_1(\tilde p))\ge v$'.
For instance, there is a sequence of compact nilpotent $n$-manifolds, $N/\Gamma_i$, which supports
no flat metric, satisfying $|\text{sec}_{N/\Gamma_i}|\le \epsilon_i\to 0$, $\op{diam}(N/\Gamma_i)=1$
and for all $\tilde x_i\in N$, $\frac{\volume(B_1(\tilde x_i))}{\svolball{0}{1}}\to 1$ uniformly (cf. \cite{Gr}).

We now state our quantitative version for Theorem 0.3.

\begin{thmD}
	Given $n, d>0$, there exists $\epsilon(n,d)>0$ such that for any $0\le \epsilon<\epsilon(n,d)$,
	if a compact $n$-manifold $M$ satisfies
	$$\op{Ric}_M\ge -(n-1), \quad d\ge \op{diam}(M),\quad h(M)\ge n-1-\epsilon,$$
	then $M$ is diffeomorphic to a hyperbolic manifold by a $\Psi(\epsilon|n,d)$-isometry.
\end{thmD}

Theorem D implies Theorem 0.3. As discussed following Theorem B, Theorem D does not
hold if one removes a bound on diameter.

To explore relations between Theorem B and Theorem D, we need the
following property:

\begin{theorem}
	Let $M_i$ be a sequence of compact $n$-manifold of $\op{Ric}_{M_i}\ge -(n-1)$ such
	that $M_i\overset{GH}{\longrightarrow}M$. If $M$ is a compact Riemannian $n$-manifold, then $h(M_i)\to h(M)$ as $i\to \infty$.
\end{theorem}

Combining Theorem B, Theorem D and Theorem 0.5, we obtain the following
corollary:

\begin{corollary}
	Let $M$ be a compact $n$-manifold such that
	$$\op{Ric}_M\ge -(n-1),\quad \op{diam}(M)\le d.$$
	Then the following conditions are equivalent as $\epsilon\to 0$:
	
	\noindent (0.6.1) $M$ is diffeomorphic and $\epsilon$-close to
	a hyperbolic manifold.
	
	\noindent (0.6.2) $\frac{\op{vol}(B_1(\tilde x))}{\svolball{-1}{1}}
	\ge 1-\epsilon$, for
	any $\tilde x\in \tilde M$.
	
	\noindent (0.6.3) $h(M)\ge n-1-\epsilon$.
\end{corollary}

A few remarks are in order:

\begin{remark} 
	Theorem A generalizes Theorem 0.2;
	first, if $M$ has an almost maximal volume, then $M$ is simply connected and thus
	$M$ satisfies the conditions of Theorem A for $\rho>\pi$. Secondly, Theorem A applies
	to all spherical $n$-space form; all but finitely many are collapsed when $n$ is odd.
	Theorem A also improves Theorem 0.2;
	if $M$ in Theorem A is simply connected, then $M$ is diffeomorphic
	and $\Psi(\epsilon|n)$-close to $S^n_1$, while the conditions do not
	apriorily imply that the volume of $M$ almost equals to $\op{vol}(S^n_1)$. We point it out
	that the case in Theorem A for $\rho>\pi$ also recovers Theorem 4 in \cite{Au}
	which is a generalization of Theorem 0.2.
\end{remark}

\begin{remark} 
	If $M$ satisfies the condition in Theorem B or Theorem D, then
	$\op{vol}(M)$ is not less than the volume of the hyperbolic metric on $M$ (\cite{BCG}),
	which is bounded below by a constant $v(n)$ (Heintze-Margulis, cf. \cite{He}).
	In particular, this answers a question in \cite{LW2} whether $M$ of almost maximal
	volume entropy can collapse.
\end{remark}

\begin{remark} 
	The gap phenomena in Theorem C that
	``$\frac{\volume(B_\rho(x^*))}{\svolball{0}{\rho}}\ge 1-\epsilon$'' implies
	that ``$\frac{\volume(B_\rho(x^*))}{\svolball{0}{\rho}}=1$ '' is related to
	the bounded ratio of diameters on $M$ and $\tilde M$ when $\pi_1(M)$ is
	finite (\cite{KW}). Nevertheless, this volume gap phenomena seems not to be explored
	before; compare with flat manifolds rigidity
	under non-negative Ricci condition (e.g., Corollary 27 and 29, \cite{Pet}).
\end{remark}

\begin{remark} 
	Note that in Theorem 0.4, $\Psi(\delta|n,\rho,v)$
	is independent of $\epsilon$; this is because a limit space of a sequence
	manifolds in Theorem 0.4 with $\delta_i\to 0$ is isometric to a flat manifold
	(see Lemma 3.8). The independence of $\epsilon$ was pointed out to us
	by S. Honda after the first version was put on ArXiv.
\end{remark}

\begin{remark} 
	Let $M$ be a compact hyperbolic $n$-manifold. The minimal
	volume rigidity in \cite{BCG} says that any metric $g$ on $M$
	of $\op{Ric}_g\ge -(n-1)$ satisfies that $\op{vol}(M,g)\ge\op{vol}(M)$,
	and ``$=$'' if and only if $g$ is the hyperbolic metric on $M$.
	By Theorem 0.3, $h(M,g)\le h(M)$ and ``$=$'' if and only if $g$ coincides
	with the hyperbolic metric. In comparing the quantitative minimal
	volume rigidity (Theorem 1.3 in \cite{BBCG}) with Theorem D, a substantial
	difference is that the former requires a non-collapsing condition
	but no condition on diameter, while the latter requires a bound on
	diameter but no non-collapsing condition.
\end{remark}

\begin{remark}
	For a special case of Theorem D that manifolds have
	strictly negative sectional curvature, see \cite{LW2}.
\end{remark}

\begin{remark} 
	If, in Theorem A, B and D, the curvature condition is replaced
	by $\op{Ric}_M\ge (n-1)H$ ($H>0$ or $H<0$), then conclusions hold with respect to
	the space form of constant curvature $H$, provided that $\epsilon$ also depends on $H$.
\end{remark}

\begin{remark}
	In the proof of Theorem A-C,  we show that the Riemannian
	universal covering space satisfies that for any $\tilde x\in \tilde M$, $\frac{\op{vol}
		(B_{\rho'}(\tilde x))}{\svolball{H}{\rho'}}
	\ge 1-\Psi(\epsilon|n,\rho,d,v)$ ($H=1, -1,$ or $0$), where $\rho'=\rho'(n,\rho,d,v)>0$, see Corollary 3.3.
\end{remark}

In the light of Theorem A-C,  we propose the following:

\begin{conjecture}[Quantitative maximal local rewinding volume rigidity]
	Given $n, \rho>0$ and $H=\pm 1$ or $0$, there exists a constant
	$\epsilon(n,\rho)>0$ such that for any $0<\epsilon<\epsilon(n,\rho)$,
	if a compact $n$-manifold $M$ satisfies
	$$\op{Ric}_M\ge (n-1)H,\quad \frac{\volume
		(B_\rho(x^*))}{\svolball{H}{\rho}}\ge 1-\epsilon,\quad\forall\, x\in M,$$
	then $M$ is diffeomorphic and $\Psi(\epsilon|n,\rho)$-close
	to a space form of constant curvature $H$, provided that $\op{diam}(M)\le d$ (and
	thus $\epsilon(n,\rho,d)$) when $H\ne 1$.
\end{conjecture}

The following is a supporting evidence for Conjecture 0.15 (see \cite{CRX}).

\begin{thmE}
	Conjecture 0.15 holds for the class of Einstein manifolds.
\end{thmE}

We now briefly describe our approach to Theorem A-C and Theorem D which is
quite involved with tools from several fields. The most significant tool is from
the Cheeger-Colding theory (\cite{Ch}, \cite{Co2}, \cite{CC1,CC2,CC3}) and the Perel'man's
pseudo-locality of Ricci flows (\cite{BW}, \cite{Ha1,Ha2}, \cite{Pe2}).
In our proof of Theorem A, we established a $C^0$-convergence (see Theorem 2.7),
and in the proof of Theorem D, we
establish that an almost volume annulus of fixed width and radius going to $\infty$ ($H\le 0$) contains a
large ball that is almost metric warped product (see Theorem 1.4).
This result complements the Cheeger-Colding's theorem that an almost
volume annulus (of bounded radius) is an almost metric annulus, and also
yields a new proof of Theorem 0.3 (see Remark 4.5) that does not rely on \cite{LiW} (cf. \cite{LW1}, \cite{Li}).

Starting with a contradicting sequence to Theorem A-C, $M_i\overset{GH}{\longrightarrow}X$, such that
$\frac{\volume(B_\rho(x_i^*))}{\svolball{H}{\rho}}\ge 1-\epsilon_i$
for all $x_i\in M_i$, and we will study the associate equivariant sequence
of the Riemannian universal
covering spaces, which satisfies the following commutative diagram (\cite{FY1}):
$$\begin{CD} (\tilde M_i,\tilde p_i,\Gamma_i)@>GH>>(\tilde X,\tilde p,G)
\\ @VV \pi_i V   @VV \pi  V \\
(M_i,p_i)@>GH>> (X,p),\end{CD} \eqno (0.16)$$
where $\Gamma_i=\pi_1(M_i,p_i)$ is the fundamental group, $G$ is the limiting Lie group (\cite{CC3})
and the identity component $G_0$ is nilpotent (\cite{KW}). We will first show that
$\tilde X$ is locally isometric to a space form. For any $\tilde x\in
\tilde X$, let $\tilde x_i\in \tilde M_i$ such that $\tilde x_i\to \tilde x$,
we study a local version of (0.16):
$$\begin{CD} (U^*_\rho,x_i^*,\Lambda_i)@
>GH>>(\tilde Y,x^*,K)\\ @VV  \tilde \pi_i^* V   @VV \tilde \pi^*  V \\
(\pi_i^{-1}(B_\rho(x_i)),\tilde x_i)@>GH>> (Y,\tilde x),\end{CD} \eqno (0.17)$$
where $\Lambda_i=\pi_1(\pi_i^{-1}(B_\rho(x_i)),\tilde x_i)$.
According to the Cheeger-Colding's theorem that an almost volume annulus
is an almost metric annulus, $\frac{\volume
	(B_\rho(x_i^*))}{\svolball{H}{\rho}}\ge 1-\epsilon_i$
implies that $d_{GH}(B_{\frac \rho2}(x_i^*), \sball{H}{\frac \rho2})<
\Psi(\epsilon_i|n,\rho)$ (Theorem 1.2), and thus $\tilde Y$ is
locally isometric to a $H$-space form. Since $\tilde M_i$ is not collapsed,
$K$ is discrete. It remains to check that $K$
acts freely (Theorem 2.1), thus a small ball at $\tilde x$ is isometric to a small ball
in an $n$-dimensional $H$-space form. If $e\ne \gamma\in K$ and $q^*\in B_{\frac \rho4}(x^*)$ such
that $\gamma(q^*)=q^*$,  under the non-collapsing equivariant convergence
we show that $\gamma$ and $q^*$ can be
chosen so that there are $\gamma_i\in \Lambda_i$ of order equal to that of $\gamma$,
$\gamma_i\to \gamma$, $q_i^*\to q^*$ and the displacement of
$\gamma_i$ at $q_i^*$, $\mu_i\to 0$,  is almost minimum around $q_i^*$. In our circumstance,
the rescaling sequence, $$\begin{CD}(\mu_i^{-1}U^*_\rho,q_i^*,\left<\gamma_i
\right>)@>GH>>(\Bbb R^n,\tilde q,\left<\gamma'\right>)\end{CD},$$ which leads
to a contradiction because $\gamma'$ must fix some point in $\Bbb R^n$, while
$\gamma_i$ moves every point at least a definite amount, where $\left<\gamma_i\right>$ denotes the subgroup
generated by $\gamma_i$.

If $G$ is discrete, similar to the above we conclude that $G$ acts freely on
$\tilde X$ (Theorem 2.1), and thus $X$ is isometric to an $n$-dimensional $H$-space form.
We then get a contradiction by applying the diffeomorphic stability
theorem in \cite{CC2}. For $H=-1$, we will show that $G$ is discrete (Theorem 2.5): using the
nilpotency of $G_0$ and the compactness of
$\tilde X/G$ we show that $G_0$ contains neither elliptic nor hyperbolic
elements (Lemma 2.6). Using (0.16), we construct a geodesic segment
in some $G_0$-orbit, and thus conclude that $G_0$ contains no parabolic element
i.e., $G_0=e$. This finishes the proof of Theorem B.

For $H=0$, $\tilde X=\Bbb R^k\times F$ and
$\tilde M_i=\Bbb R^k\times N_i$ (Cheeger-Gromoll splitting theorem),
where $F$ is a compact flat manifold, and $N_i$ is a compact simply connected
manifold of non-negative Ricci curvature.
We show that $\op{diam}(N_i)$ is uniformly bounded above, and thus
applying the diffeomorphic stability theorem in \cite{CC2} we derive a
contradiction.

For $H=1$, in (0.16)  we may assume an $\epsilon_i$-equivariant diffeomorphism,
$\tilde h_i: (\tilde M_i,\Gamma_i)\to (S^n_1,G)$ (\cite{CC2}). Via $\tilde h_i$, we identify
$(M_i,\Gamma_i)$ as a free $\Gamma_i$-action on $S^n_1$ by $\epsilon_i$-isometries.
By \cite{MRW},  for $i$ large there is an injective homomorphism, $\phi_i: \Gamma_i\to G$ (see
Lemma 3.4).  We show that the $\phi_i(\Gamma_i)$-action on $S^n_1$ is free (see (3.5.1)).
By now we can perform the center of mass to perturb $\op{id}_{S^n_1}$ to
a map, $\tilde f_i: S^n\to S^n$, that commutes the $\Gamma_i$-action
with the $\phi_i(\Gamma_i)$-action. It remains to show that $\tilde f_i$ is a diffeomorphism,
and thus a contradiction. According to \cite{GK},  $\tilde f_i$ is a diffeomorphism when
the $\Gamma_i$- and $\phi_i(\Gamma_i)$-actions are close in $C^1$-norm. To see it,
we will
use Ricci flows of $\tilde g_i$: using Perel'man's pseudo-locality
(\cite{Pe2}) and a distance estimate in \cite{BW} we show that a solution
$\tilde g_i(t)$ is $C^0$-close to $\b g^1$ on $S^n_1$ (see Theorem 2.7); which is also
locally $C^{1,\alpha}$-close to $\b g^1$ up to a definite rescaling.
Since $\Gamma_i$ remains to be isometries with respect to $\tilde g_i(t)$,
the above regularities guarantee the desired $C^1$-closeness (see (3.5.2)).

In the proof of Theorem D, we again start with a contradicting sequence as in (0.16),
and it suffices to show that $\tilde X$ is isometric to
$\Bbb H^n$, and by the volume convergence (\cite{Co2}) $M_i$ satisfies the conditions
of Theorem B, a contradiction. Fixing $R>50d$, we will
prove that
$d_{GH}(B_R(\tilde p_i), \sball{-1}{R})<\Psi (\epsilon_i|n,d,R)$,  where $\sball{-1}{R}$ is
a ball in $\Bbb H^k$ for some $k\le n$
(Lemma 4.4). First, following \cite{Li} we show that $h(M)\ge n-1-\epsilon$
implies a sequence, $r_i\to \infty$, such that the ratio, $\lim_{i\to \infty}\frac {\op{vol}
	(\partial B_{r_i+50R}(\tilde p))}{\op{vol}(\partial B_{r_i-50R}(\tilde p))}
\ge e^{100R(n-1-\epsilon)}$, which approximates the limit of the same type ratio on
$\Bbb H^n$. Because $\op{vol}(A_{r_i-50R,r_i+50R}(\tilde p))\to \infty$ as $r_i\to \infty$,
the Cheeger-Colding's theorem that an almost volume annulus is an almost metric
annulus cannot be applied in our situation. Instead, we establish the following (weak) property
(see Theorem 1.4):
annulus $A_{r_i-50R,r_i+50R}(\tilde p)$ contains a ball $B_{2R}(\tilde q_i)$ such
that $d_{GH}(B_{2R}(\tilde q_i), \sball{-1}{2R})<\Psi (\epsilon_i,r_i^{-1}|n,R)$, which
leads to the desired estimate via pullback $B_{2R}(\tilde q_i)$ to $B_{2R}(\gamma_i(\tilde q_i))
\supseteq B_R(\tilde p)$ with suitable element $\gamma_i\in \Gamma_i$.

The remaining proof is to show that $k=n$.  If $k<n$,  then $M_i$ is collapsed.
By \cite{FY1}  and \cite{FY2} (see Lemma 1.13), there is $\epsilon>0$ such that the subgroup
$\Gamma_i^\epsilon\subset \Gamma_i$ generated by elements whose displacement on
$B_1(\tilde p_i)$ are uniformly smaller than $\epsilon$ converges to $G_0$.
From the proof of Theorem B, $G_0$ is trivial and thus $\Gamma^\epsilon_i$ is finite.
Since $h(M_i)$ can be calculated in terms of the growth of $\pi_1(M_i)$
at $\tilde p_i$, via center of mass method we construct an almost  $\Gamma_i/
\Gamma_i^\epsilon$-conjugate map from $(\tilde M_i/\Gamma_i^\epsilon,
\Gamma_i/\Gamma_i^\epsilon)\to
(\Bbb H^k,G)$ which is also an $\epsilon_i$-Gromov-Hausdorff approximation
when restricting to $B_R(\tilde p_i)$ (Lemma 4.7), we are able to estimate
$h(M_i)\le k-1+\epsilon_i$ (Theorem 4.6), a contradiction.

The rest of the paper is organized as follows:

In Section 1, we supply basic notions and tools concerning a convergent sequence
of compact $n$-manifolds with Ricci curvature bounded below and diameter bounded above,
which will be freely used through the rest of the paper. In particular, we will state
our result that an asymptotic volume annulus contains many disjoint balls of almost
warped product structure (see Theorem 1.4), which provides information complements to
the Cheeger-Colding's theorem that almost volume annulus is almost metric annulus (Theorem 1.3).

In section 2, we will establish three key properties for our proofs of Theorems A-C and D:
a sufficient condition for a limiting group $G$ to act freely on a limit
space $\tilde X$ (Theorem 2.1), for $H=-1$, $G$ is discrete (Theorem 2.5) and a
$C^0$-convergence of Ricci flows associate to a sequence of GH-convergence
with Ricci curvature bounded below (Theorem 2.7).

In Section 3, we will prove Theorem A-C, Theorem E and Theorem 0.4.

In Section 4, we will prove Theorem D by assuming Theorem 1.4.
We will also prove Theorem 0.5 and Corollary 0.6.

In Section 5, we will prove Theorem 1.4.

The authors would like to thank Binglong Chen for a helpful discussion
on Ricci flows.

\section{Preliminaries} 
The purpose of this section is to supply notions and basic properties from
the fundamental work of Cheeger-Colding on degeneration of Riemannian metrics with
Ricci curvature bounded from below, as well as those related to equivariant
Gromov-Hausdorff convergence. These will be used through out this paper, and
we refer the readers to \cite{Ch}, \cite{CC1,CC2,CC3}, \cite{Co1,Co2}, \cite{FY1, FY2} for details.

We will also state our result that an almost volume annulus of fixed width and
large radius contains many disjoint balls with almost warped product structure
(see Theorem 1.4).

\subsection*{a. Manifolds of Ricci curvature bounded below}
~

Let $N$ be a Riemannian $(n-1)$-manifold, let $k:(a,b)\to \Bbb R$ be a smooth positive
function and let $(a,b)\times_k N$ be the $k$-warped product whose Riemannian tensor is
$$g = dr^2 + k^2(r)g_N.$$
The Riemannian distance $|(r_1,x_1)(r_2,x_2)|$ ($x_1\neq x_2$) equals to the infimum of the length
$$\int_0^l\sqrt{ (c_1'(t))^2+k^2(c_1(t))}dt$$
for any smooth curve $c(t)=(c_1(t),c_2(t))$ such that $c(0)=(r_1,x_1)$, $c(l)=(r_2,x_2)$
and $|c_2'|\equiv 1$, and $|(r_1,x)(r_2,x)|=|r_2-r_1|$.
Thus given $a,b,k$, there is a function (e.g., the law of cosine on space forms)
$$\rho_{a,b,k}(r_1,r_2,|x_1x_2|)=|(r_1,x_2)(r_2,x_2)|.$$
Using the same formula for $|(r_1,x_1)(r_2,x_2)|$,
one can extend the $k$-warped product $(a,b)\times_k Y$ to any metric space $Y$ (not necessarily
a length space); see \cite{CC1}.

We first recall the following Cheeger-Colding's ``almost volume
warped product implies
almost metric warped product'' theorem.

\begin{theorem}[\cite{CC1}]
Let $M$ be a Riemannian manifold, let $r$ be a distance
function to a compact subset in $M$, let $0<\alpha'<\alpha, \alpha-\alpha'>2\xi>0$, let $A_{a,b}=r^{-1}((a,b))$ and let
$$\Cal V(\xi)=\inf\left\{\left.\frac{\volume (B_{\xi}(q))}{\volume (A_{a,b})}
\right| \text{ for all $q\in A_{a,b}$ with $B_{\xi}(q)\subset A_{a,b}$ }\right\}.$$
If
$$\op{Ric}_M \ge -(n-1)\frac{k''(a)}{k(a)}\qquad (\text{on $r^{-1}(a)$}),$$
$$\Delta r\le (n-1)\frac{k'(a)}{k(a)} \qquad (\text{on $r^{-1}(a)$}),$$
$$\frac{\volume(A_{a,b})}{\volume(r^{-1}(a))}\ge (1-\epsilon)\frac{\int_a^bk^{n-1}(r)dr}{k^{n-1}(a)}. \eqno (1.1.1)$$
Then there exists a length metric space $Y$, with at most $\#(a,b,k,\Cal{V}(\xi))$
components $Y_i$, satisfying
$$\op{diam}(Y_i)\le D(a,b,k,\Cal V(\xi)),$$
such that $$d_{GH}(A_{a+\alpha,b-\alpha},(a+\alpha,b-\alpha)\times_kY)\le
\Psi(\epsilon|n,k,a,b,\alpha',\xi,\Cal V(\xi)) \eqno (1.1.2)$$
with respect to the two metrics $d^{\alpha',\alpha}$ and $\b d^{\alpha',\alpha}$, where
$d^{\alpha',\alpha}$ (resp. $\b d^{\alpha',\alpha}$) denotes the restriction of the
intrinsic metric of $A_{a+\alpha',b-\alpha'}$ on $A_{a+\alpha,b-\alpha}$ (resp.
$(a+\alpha',b-\alpha')\times_kY$ on $(a+\alpha,b-\alpha)\times_kY$).	
\end{theorem}

Let $$\op{sn}_H(r)=\begin{cases} \frac{\sin \sqrt Hr}{\sqrt H} & H>0\\ r & H=0\\ \frac{\sinh \sqrt{-H}r}{\sqrt{-H}} & H<0\end{cases}.$$
Applying Theorem 1.1 to $\op{sn}_H(r)$ with $r(x)=d(p,x): M\to \Bbb R$,
we conclude the following ``almost maximal volume ball implies almost
space form ball'', which is important to our work (one may need to shift the center a bit to
see the following).

\begin{theorem}
For $n, \rho, \epsilon >0$, if a complete $n$-manifold $M$ contains a point $p$ satisfies
$$\op{Ric}_M\ge (n-1)H,\quad \frac{\volume(B_\rho(p))}{\svolball{H}{\rho}}\ge 1-\epsilon,$$
then $d_{GH}(B_{\frac \rho2}(p),\sball{H}{\frac \rho2})<\Psi(\epsilon|n,\rho,H)$.	
\end{theorem}

Another important application of Theorem 1.1 is the following an ``almost
volume annulus'' is an ``almost metric annulus''. For $p\in M$, $L>2R>0$,
let $A_{L-2R,L+2R}(p)=\{x\in M, \,\,\, L-2R<|xp|<L+2R\}$.

\begin{theorem}
Given $n, H\leq 0, L>2R>0$, if a complete $n$-manifold
$M$ contains a point $p$ satisfies
$$\op{Ric}_M\geq (n-1)H,\quad \frac{\volume(\partial B_{L-2R}(p))}{\svolsp{H}{{L-2R}}}
\leq (1+\epsilon)\frac{\volume( A_{L-2R, L+2R}(p))}{\svolann{H}{L-2R, L+2R}}, \eqno (1.3.1)$$
then
$$d_{GH}(A_{L-R,L+R}(p), (L-R,L+R)\times_{\op{sn}_H(r)}Y)\leq \Psi(\epsilon| n, L, R, H), \eqno(1.3.2)$$
where $Y$ is a length metric space (may be not connected).
\end{theorem}

It turns out that in our proof of Theorem D, the condition that $h(M)
\ge n-1-\epsilon$ implies that (1.3.1) is satisfied asymptotically i.e.,
only as $L\to \infty$ (see Lemma 4.2). Because in our circumstance
$\op{vol}(A_{L-R,L+R}(p))$ $\to \infty$ as $L\to \infty$, it is not
possible to have (1.3.2) in our circumstance.

In our proof Theorem D, it is crucial for us to establish the following result.

\begin{theorem}
Given $n, H\leq 0,  L >> R\ge \rho> 0, \epsilon > 0$, there exists a constant
$c = c(n,H,R,\rho)$ such that if a complete $n$-manifold $M$ contains a point $p$ satisfies
(1.3.1), then there are disjoint $\rho$-balls, $B_{\rho}(q_i)\subset A_{L-R, L+R}(p)$, for each $B_{\rho}(q_i)$,
$$ d_{GH}(B_{\rho}(q_i)), B_{\rho}((0,x_i))\leq \Psi(\epsilon, L^{-1}| n, H, R, \rho) \eqno (1.4.1) $$
where $B_{\rho}((0,x_i))\subset \Bbb{R}^1\times_{e^{\sqrt{-H}r}} Y_i$ for some length metric space $Y_i$, and
$$ \quad \frac{\volume(\bigcup_i B_{\rho}(q_i))}{\volume(A_{L-R, L+R}(p))}\geq c(n,H,R,\rho). \eqno (1.4.2)$$
In particular, for $H=0$, we have that each $B_{\rho}(q_i)$ is almost splitting.
\end{theorem}

Roughly, Theorem 1.4 says that for any fixed $R>0$, if $A_{L-2R,L+2R}(p)$ is
an almost volume annulus as $L\to \infty$, then (even if its volume blows up
to infinity) one can have lots of disjoint balls of fixed radius
$\rho\le R$ in the annulus, each of which is close to a ball
in a metric annulus.

The proof of Theorem 1.4 uses the same techniques from \cite{Ch} and \cite{CC1}, and
because it is technical and tedious, we will leave the proof in section 5.

\begin{remark}
The almost volume annulus condition (1.3.1) implies the
following:
$$\frac{\volume(\partial B_{L+R}(p))}{\svolsp{H}{{L+R}}}\geq (1-\Psi(\epsilon| n, H, R))
\frac{\volume(\partial B_{L-R}(p))}{\svolsp{H}{L-R}} \eqno (1.5.1)$$
From the proof of Theorem 1.3 in \cite{CC1}, one sees that indeed only (1.5.1) is applied. Furthermore,
(1.3.1) and (1.5.1) are equivalent conditions when $\epsilon$ is small.
\end{remark} 

Consider a sequence of complete $n$-manifolds, $(M_i,p_i)\overset{GH}{\longrightarrow}(X,p)$, such that
$\op{Ric}_{M_i}\ge -(n-1)$. If $M_i$ is not collapsed, then a basic property is:

\begin{theorem}[\cite{Co2,CC2}]
Let $(M_i,p_i)\overset{GH}{\longrightarrow}(X,p)$ such that $\op{Ric}_{M_i}$
$\ge -(n-1)$. If $\op{vol}(B_1(p_i))\ge v>0$, then for any $r>0$, $M_i\ni x_i\to x\in X$,
$\op{vol}(B_r(x_i))\to \op{Haus}^n(B_r(x))$, where $\op{Haus}^n$ denotes the
$n$-dimensional Hausdorff measure.
\end{theorem}

Let $X$ be a complete separable length metric space. A
point $x\in X$ is called a $(\epsilon,r)$-Reifenberg point, if for any $0<s<r$,
$$d_{GH}(B_s(x),\sball{0}{s})\le \epsilon s.$$
$X$ is called a $(\epsilon,r)$-Reifenberg space if every point in $X$ is a
$(\epsilon,r)$-Reifenberg point.

\begin{theorem}[\cite{CC2}]
Let $M_i\overset{GH}{\longrightarrow}X$ be a sequence of complete $n$-manifolds
of $\text{Ric}_{M_i}\ge -(n-1)$, and $X$ is compact. Then there is a constant $\epsilon(n)>0$ such that
for $i$ large

\noindent (1.7.1) If $X$ is a  $(\epsilon,r)$-Riefenberg space with $\epsilon<\epsilon(n)$,
then there is a homeomorphic bi-H\"older equivalence between $M_i$ and $X$.

\noindent (1.7.2) If $X$ is a Riemannian manifold, then there is a diffeomorphic bi-H\"older
equivalence between $M_i$ and $X$.	
\end{theorem}

\begin{theorem}[\cite{CC3}]
Let $(M_i,p_i)\overset{GH}{\longrightarrow}(X,p)$ such that $\op{Ric}_{M_i}
\ge -(n-1)$. If $\op{vol}(B_1(p_i))\ge v>0$, then the isometry group of
$X$ is a Lie group.
\end{theorem}

Theorem 1.8 holds for any limit space of Riemannian $n$-manifolds with Ricci curvature bounded below (\cite{CN}).

According to the classical Margulis Lemma, if $M$ is a symmetric space, the subgroup
of the fundamental group of $M$ generated by loops of small length is virtually
nilpotent. Magulis Lemma was extended in \cite{FY1} to manifolds of $\op{sec}\ge -1$
that the subgroup is virtually nilpotent,  and in \cite{KPT} a bound on the index of the
nilpotent subgroup was obtained depending only on $n$.  Recently, Kapovitch-Wilking
proved the following generalized Magulis Lemma (conjectured by Gromov):

\begin{theorem}[\cite{KW}]
There are constants $\epsilon(n), w(n)>0$ if
$M$ is a complete $n$-manifold of $\op{Ric}_M\ge -(n-1)$, $p\in M$, then the
image subgroup, $\op{Im}(\pi_1(B_\epsilon(p))\to \pi_1(M))$
contains a nilpotent subgroup of index $\le w(n)$, with the nilpotent basis of length at most $n$.
\end{theorem}

\subsection*{b. Equivariant Gromov-Hausdorff convergence}
~

The reference of this part is \cite{FY1}, \cite{FY2} , \cite{KW} (cf. \cite{Ro2}).

Let $X_i\overset{GH}{\longrightarrow}X$ be a convergent sequence of compact length metric spaces, i.e.,
there are a sequence $\epsilon_i\to 0$ and a sequence of
maps $h_i: X_i\to X$, such that $||h_i(x_i)h_i(x_i')|_X-|x_ix_i'|_{X_i}|
<\epsilon_i$ ($\epsilon_i$-isometry), and for any $x\in X$, there is $x_i\in X_i$ such that
$|h_i(x_i)x|_X<\epsilon_i$ ($\epsilon_i$-onto), and $h_i$ is called an $\epsilon_i$-Gromov-Hausdorff
approximation, briefly, $\epsilon_i$-GHA. From now on, we will omit the subindex in the distance
function ``$|\cdot\cdot|$''.

Assume that $X_i$ admits a closed group $\Gamma_i$-action by isometries. Then
$(X_i,\Gamma_i)\overset{GH}{\longrightarrow}(X,\Gamma)$ means that there are a sequence $\epsilon_i\to 0$ and
a sequence of $(h_i,\phi_i,\psi_i)$, $h_i: X_i\to X$, $\phi_i:
\Gamma_i\to \Gamma$ and $\psi_i: \Gamma\to \Gamma_i$ which are $\epsilon_i$-GHAs such that for all
$x_i\in X_i, \gamma_i\in \Gamma_i$ and $\gamma\in \Gamma$,
$$\begin{gathered}
|h_i(x_i)[\phi_i(\gamma_i)h_i(\gamma_i^{-1}(x_i))]|<\epsilon_i,\\
|h_i(x_i)[\gamma^{-1}(h_i(\psi_i(\gamma)(x_i)))]|<\epsilon_i,
\end{gathered} \eqno (1.10)$$
where $\Gamma$ is a closed group of isometries on $X$, $\Gamma_i$ and
$\Gamma$ are equipped with the induced metrics from $X_i$ and $X$.
We call $(h_i,\phi_i,\psi_i)$
an $\epsilon_i$-equivariant GHA.

When $X$ is not compact, then the above notion of equivariant convergence naturally extends to a pointed version $(h_i,\phi_i,\psi_i)$:
$h_i: B_{\epsilon_i^{-1}}(p_i)\to B_{\epsilon_i^{-1}+\epsilon_i}(p)$, $h_i(p_i)=p$,
$\phi_i: \Gamma_i(\epsilon_i^{-1})\to \Gamma(\epsilon_i^{-1}+\epsilon_i)$, $\phi_i(e_i)=e$,
$\psi_i: \Gamma(\epsilon_i^{-1})\to \Gamma_i(\epsilon_i^{-1}+\epsilon_i)$, $\psi_i(e)=e_i$,
and (1.10) holds whenever the multiplications stay in the domain of $h_i$, where
$\Gamma_i(R)=\{\gamma_i\in\Gamma_i,\,\,|p_i\gamma_i(p_i)|\le R\}$.

\addtocounter{theorem}{1} 
\begin{lemma}
Let $(X_i,p_i)\overset{GH}{\longrightarrow}(X,p)$, where $X_i$ is a  complete
locally compact length space. Assume that $\Gamma_i$ is a closed group of isometries
on $X_i$. Then there is a closed group $G$ of isometries on $X$
such that passing to a subsequence, $(X_i,p_i,\Gamma_i)\overset{GH}{\longrightarrow}(X,p,G)$.
\end{lemma}

\begin{lemma}
Let $(X_i, p_, \Gamma_i)\overset{GH}{\longrightarrow}(X,p,G)$, where $X_i$ is a complete
locally compact length space and $\Gamma_i$ is a closed subgroup of isometries.
Then $(X_i/\Gamma_i,\bar p_i)\overset{GH}{\longrightarrow}(X/G,\bar p)$.
\end{lemma}

For $p_i\in X_i$, let $\Gamma_i=\pi_1(X_i,p_i)$ be the fundamental group.
Assume that the universal covering space, $\pi_i: (\tilde X_i, \tilde p_i)\to (X_i,p_i)$,
exists.

\begin{lemma}
Let $X_i\overset{GH}{\longrightarrow}X$ be a sequence of compact length metric
space. Then passing to a subsequence such that the following
diagram commutes,
$$\begin{CD} (\tilde X_i,\tilde p_i,\Gamma_i)@>GH>>(\tilde X,\tilde p,G)
\\ @VV \pi_i V   @VV \pi  V \\
(X_i,p_i)@>GH>> (X,p).\end{CD}$$
If $X$ is compact and $G/G_0$ is discrete, then there is $\epsilon>0$ such that the subgroup,
$\Gamma_i^\epsilon$, generated by elements with displacement bounded
above by $\epsilon$ on $B_{2d}(\tilde p_i)$, is normal and for $i$ large,
$\Gamma_i/\Gamma_i^\epsilon\overset{\op{isom}}{\cong} G/G_0$.
\end{lemma}

Combining Lemma 1.12 and 1.13, we obtain the following commutative diagram:
$$\begin{CD} (\tilde X_i,\tilde p_i,\Gamma_i)@>GH>>(\tilde X,\tilde p,G)
\\ @VV \hat \pi_i V   @VV \hat \pi  V \\ (\hat X_i,\hat p_i,\hat \Gamma_i)@>GH>>(\hat X,\hat p,\hat G)
\\ @VV \bar \pi_i V   @VV \bar \pi  V \\
(X_i,p_i)@>GH>> (X,p),\end{CD} \eqno (1.14)$$
where $\hat X_i=\tilde X_i/\Gamma_i^\epsilon, \hat X=\tilde X/G_0,\hat \Gamma_i=\Gamma_i/\Gamma_i^\epsilon$
and $\hat G=G/G_0$.


\section{The Free Action, The Discreteness of Limiting Groups and The $C^0$-convergence}

In this section, we will establish three key properties for our proofs of
Theorems A-D: Theorem 2.1, Theorem 2.5 and Theorem 2.7.

\subsection*{a. Free limit isometric actions}
~

Let $(M_i,p_i)$ be a sequence of complete $n$-manifolds, let $\pi_i^*: (U^*_d,p_i^*)\to (B_d(p_i),p_i)$ be the Riemannian universal  covering spaces,
and let $\Lambda_i=\pi_1(B_d(p_i),p_i)$ denote the fundamental group.

\begin{theorem}
Given $n, d, v, r>0$, there exists a constant $\epsilon=\epsilon(n,v)>0$
such that if a sequence of complete $n$-manifolds, $(M_i,p_i)$, satisfies
$$\text{Ric}_{M_i}\ge -(n-1),\,\op{vol}(B_1(p_i))\ge v,\, \forall\, x^*\in
B_{\frac{d}{2}}(p^*)\text{  is a $(\epsilon,r)$-Reifenberg point}$$
and the following commutative diagram:
$$\begin{CD} (U_d^*,p_i^*,\Lambda_i)@>GH>>(\tilde X,p^*,K)
\\ @VV \pi_i^* V   @VV \pi^*  V \\
(B_d(p_i),p_i)@>GH>> (B_d(p),p),\end{CD}$$
then the discrete group $K$ acts freely on $B_{\frac d4}(p^*)$ i.e., $K$ has no
isotropy group in $B_{\frac d4}(p^*)$.
\end{theorem}

\begin{corollary}
Given $n, \rho, v>0$ and $H\ge -1$, there exists
a constant $\epsilon=\epsilon(n,v)>0$ such that if a sequence of
complete $n$-manifolds, $(M_i,p_i)$, satisfies
$$\op{Ric}_{M_i}\ge (n-1)H_i\to (n-1)H,\quad \op{vol}(B_1(p))\ge v,
\quad\frac{\volume(B_\rho( p_i^*))}{\svolball{H}{\rho}}\ge 1-\epsilon,$$
and the following commutative diagram:
$$\begin{CD} (U_\rho^*,p_i^*,\Lambda_i)@>GH>>(\tilde X,p^*,K)
\\ @VV \pi_i^* V   @VV \pi^*  V \\
(B_\rho(p_i),p_i)@>GH>> (B_\rho(p),p),\end{CD} \eqno (2.2.1)$$
where $\pi_i^*: (U_\rho^*,p_i^*)\to (B_\rho(p_i),p_i)$ is the Riemannian universal
cover, and $\Lambda_i=\pi_1(B_\rho(x_i),p_i)$.
Then $K$ acts freely on $B_{\frac \rho4}(p^*)$ i.e., $K$
has no isotropy group in $B_{\frac \rho4}(p^*)$.
\end{corollary}

In the proof, we will use  the following lemma due to \cite{PR}:

\begin{lemma}
Let $(M_i,p_i)\overset{GH}{\longrightarrow}(X,p)$ be a sequence of complete
$n$-manifolds satisfying
$$\op{Ric}_{M_i}\ge (n-1)H_i\to (n-1)H,\quad \volume(B_\rho(p_i^*))\ge v>0,$$
and the commutative diagram (2.2.1).
If a subgroup $K_i$ of $\Lambda_i$ satisfies that $K_i\to e\in K$,
then for $i$ large, $K_i=e$.
\end{lemma}

\begin{proof}
Arguing by contradiction,  without loss of generality we may assume
$e\ne \gamma_i\in K_i$ for all $i$ such that the following diagram commutes:
$$\begin{CD} (U_\rho^*,p_i^*,\left<\gamma_i\right>)@>GH>>(\tilde X,p^*,e)
\\ @VV \hat \pi_i^* V   @VV \hat \pi^*  V \\
(U_\rho^*/\left<\gamma_i\right>,\hat p_i)@>GH>> (\hat X,\hat p),\end{CD}$$
where $\left<\gamma_i\right>$ denotes the subgroup generated by $\gamma_i\in \Lambda_i$.
Since $\left<\gamma_i\right>\overset{GH}{\longrightarrow}e$, by Lemma 1.12 $\tilde X=\hat X$,
$B_r(p_i^*)$ and $\gamma_i(B_r(p_i^*))$ $\subset B_{r+\epsilon_i}
(p_i^*)$ for some $\epsilon_i\to 0$. Let $D_i$ denote a (Dirichlet) fundamental
domain of $U_\rho^*(p_i)/\left<\gamma_i\right>$ at $p_i^*$. Then
for $0<r<\frac \rho2$, $[B_r(p_i^*)
\cap D_i] \cap [\gamma_i(B_r(p_i^*)\cap D_i)]=\emptyset$. Since $\volume
(B_\rho(p_i^*))\ge v>0$, we are able to apply Theorem 1.6 to derive
$$\begin{aligned}
\op{Haus}^n(B_r(p^*))&=\op{Haus}^n(B_r(\hat p))
=\lim_{i\to \infty}\volume(B_r(\hat p_i))\\&=\lim_{i\to \infty}
\volume(B_r(p_i^*)\cap D_i)\\
&=\lim_{i\to \infty}\frac 12
[\volume(B_r(p_i^*)\cap D_i)+\volume(\gamma_i(B_r(p_i^*)
\cap D_i))]\\
&\le \lim_{i\to \infty}\frac 12\volume(B_{r+\epsilon_i}
(p_i^*))=\frac 12\op{Haus}^n(B_r(p^*)),
\end{aligned}$$
a contradiction.
\end{proof}


\begin{proof}[Proof of Theorem 2.1]
~

Arguing by contradiction, assume a sequence, $(\epsilon_j,r_j)\to (0,0)$, and for each $j$, there
is a contradicting sequence $(M_{i,j},p_{i,j})$ to Theorem 2.1,
$$\begin{CD} (U_d^*,p_{i,j}^*,\Lambda_{i,j})@>GH>>(\tilde Y_j,p_j^*,K_j)
\\ @VV \pi_{i,j} V   @VV \pi_j  V \\
(B_d(p_{i,j}),p_{i,j})@>GH>> (B_d(p_j),p_j),\end{CD}$$
such that
$$\op{Ric}_{M_{i,j}}\ge -(n-1),\quad \op{vol}(B_1(p_{i,j}))\ge v,$$
any point $x_{i,j}^*\in B_{\frac d2}(p_{i,j}^*)$ is a $(\epsilon_j,r_j)$-Reifenberg point,
and $K_j$ has an isotropy group in $B_{\frac d4}(p_j^*)$.
Passing to a subsequence, we may assume
$$\begin{CD}(\tilde Y_j,p_j^*,K_j)@>GH>>(\tilde Y,p^*,K).\end{CD}$$
Assume $e_j\ne \gamma_j\in K_j$, $q_j^*\in B_{\frac d4}(p_j^*)$ such that $\left<\gamma_j
\right>(q_j^*)=q_j^*$. Passing to a subsequence, we may assume
$\left<\gamma_j\right>\to W$ and $q_j^*\to q^*$ such that
$W(q^*)=q^*$. We observe that Lemma 2.3 can still apply to
the above sequence i.e., if $\gamma_j\in K_j$ such that $\left<\gamma_j\right>\to e$,
then $\gamma_j=e$ for $j$ large. Hence $W\ne e$.

Without loss of generality, we may assume that $q_j^*=p_j^*$. For $e\ne \gamma\in W$,
$\gamma(p^*)=p^*$. By a standard diagonal argument, we may assume a convergent subsequence,
$$\quad \begin{CD} (U_d^*,p_{i_j,j}^*,\Lambda_{i_j,j})@>GH>>(\tilde Y,p^*,K)
\\ @VV \pi_{i_j,j} V   @VV \pi  V \\
(B_d(p_{i_j,j}),p_{i_j,j})@>GH>> (B_d(p),p).\end{CD}$$
Since $\volume(B_d(p_{i_j,j}))\ge v$, $\dim(\tilde Y)=n$ and $K$ is a Lie group
(Theorem 1.8), and therefore $K$ is discrete. Since the isotropy group $K_{p^*}$
is compact, $K_{p^*}$ is finite. Since $\gamma\in W\subset K_{p^*}$, we may
assume the order $o(\gamma)=k<\infty$.

Let $\gamma_{i,j}\to \gamma_j$. Observe that for each fixed $r_j$, $\frac{|p_{i,j}^*\gamma_{i,j}(p_{i,j}^*)|}
{r_j}\to 0$ as $i\to \infty$. We may assume the above subsequence is chosen so that
$$\frac{|p_{i_j,j}^*\gamma_{i_j,j}(p_{i_j,j}^*)|}{r_j}\le j^{-1}. \eqno (2.1.1)$$
For the sake of simple notation, from now on we will use $i=j=(i_j,j)$.

Let $\gamma_i\in \Lambda_i$ such that $\gamma_i\to \gamma$. Since for all $m\in \Bbb Z$,
$\gamma_i^m\to \gamma^m\in \{\gamma, ..., \gamma^k=e\}$, and since $K$ is discrete,
we conclude that $\left<\gamma_i\right>\to \left<\gamma\right>$ and $o(\gamma_i)=k$
(otherwise, the subgroup, $\left<\gamma_i^k\right>\to e$, a contradiction to
Lemma 2.3; compare to Remark 2.4).

Observe that if the displacement function of $\gamma_i$, $d_{\gamma_i}(z_i^*)=|z_i^*\gamma_i(z_i^*)|$, achieves a minimum
at $p_i^*$, then from the limit, $(d_{\gamma_i}(p_i^*)^{-1}U^*_d(p_i^*),p_i^*,\left<\gamma_i\right>)$,
as $i\to \infty$, one easily sees a contradiction (see below). To overcome the trouble that $d_{\gamma_i}$
may take minimum near the boundary, we claim the following property:

\noindent (2.1.2) For each $i$, there is $q_i^*\in B_{200k\cdot d_{\gamma_i}(p_i^*)}
(p_i^*)$ such that $d_{\gamma_i}(q_i^*)\le d_{\gamma_i}(p_i^*)$ and any $x_i^*\in B_{100k\cdot d_{\gamma_i}(q_i^*)}
(q_i^*)$, $d_{\gamma_i}(x_i^*)\ge \frac 1{100}\cdot d_{\gamma_i}
(q_i^*)$.

Assuming (2.1.2), we will derive a contradiction as follows: Since
$q_i^*\to p^*$ and
$d_{\gamma_i}(q_i^*)\to 0$, passing to a subsequence, we may assume
$$\begin{CD}(d_{\gamma_i}(q_i^*)^{-1}U_d^*,
q_i^*,\left<\gamma_i\right>)@>
GH>>(\tilde Y',\tilde q',\left<\gamma'\right>)\end{CD}$$
such that $\op{Ric}_{d_{\gamma_i}(q_i^*)^{-1}\tilde M_i}\ge -(n-1)d_{\gamma_i}(q_i^*)^2\to 0$.
Since points in $B_{\frac d4}(p_i^*)$ are
$(\epsilon_i,r_i)$-Reifenberg points, by
(2.1.1) we can conclude that $\tilde Y'$ is isometric to $\Bbb R^n$.
Since $o(\gamma_i)=k$, $o(\gamma')=k$ and thus $\gamma'$ has a fixed point $\tilde z'$ of distance from $\tilde q'$ at most $10k$ ($\tilde z'$ may be chosen as the center of mass for $\left<\gamma'\right>(\tilde q')$).
On the other hand, the choice of $q_i^*$ with the assigned property implies that
$d_{\gamma_i}\ge \frac 1{100}$ on $B_{100k}(q_i^*)$ (after
scaling), a
contradiction.

Verification of (2.1.2): arguing by contradiction, the failure of (2.1.2)
implies that there is $(p_i^*)_1\in
B_{100k\cdot d_{\gamma_i}(p_i^*)}(p_i^*)$ such that
$d_{\gamma_i}((p_i^*)_1)<\frac 1{100}\cdot d_{\gamma_i}(
p_i^*)$. Because $(p_i^*)_1$ lies in
$B_{200k\cdot d_{\gamma_i}(p_i^*)}(p_i^*)$, there is
$(p_i^*)_2 \in B_{100k\cdot d_{\gamma_i}((p_i^*)_1)}((p_i^*)_1)$ such that $d_{\gamma_i}((p_i^*)_2)<\frac 1{100}\cdot
d_{\gamma_i}((p_i^*)_1)<\frac{1}{100^2}d_{\gamma_i}(p_i^*)$.
Repeating the process, one gets a sequence of
points $(p_i^*)_j$ such that $(p_i^*)_j \in B_{100k\cdot
	d_{\gamma_i}((p_i^*)_{j-1})}((p_i^*)_{j-1})$ and
$d_{\gamma_i}((p_i^*)_j)<\frac{1}{100^j}d_{\gamma_i}(p_i^*)$.
Since $(p_i^*)_j\in B_{200k\cdot d_{\gamma_i}(p_i^*)}
(p_i^*)$ and the displacement of $\gamma_i$ has a positive infimum
on $B_{200k\cdot d_{\gamma_i}(p_i^*)}(p_i^*)$,
this process has to end at a finite step, a contradiction.
\end{proof} 

\begin{remark}
Note that the $\volume(B_1(p_i))\ge v>0$ is equivalent
to that the limit group $K$ is discrete, which guarantees that when $\gamma_i\to \gamma$,
$o(\gamma_i)=o(\gamma)$ for $i$ large. This does not hold if $K$ is not discrete. For instance,
let $S^1_i$ be a sequence of circle subgroup of a maximal torus $T^2$ of $O(4)$
such that $\op{diam}(T^2/S^1_i)\to 0$. Let $\Bbb Z_{q_i}\subset S^1_i$ such that $\op{diam}
(S^1_i/\Bbb Z_{q_i})\to 0$, where $q_i$ is a prime number. Since $T^2$ has no fixed point on
$S^3_1$ and $\op{diam}(T^2/S^1_i)\to 0$, $S^1_i$ has not fixed point on $S^3_1$, and therefore,
$q_i$ can be chosen so that $\Bbb Z_{q_i}$ acts freely on $S^3_1$,
and $(S^3_1,\Bbb Z_{q_i})\overset{GH}{\longrightarrow}(S^3_1,T^2)$. Since $T^2$ has a circle isotropy subgroup,
we may assume $p\in S^3_1$ and $\gamma\in T^2$ of order $2$ such that $\gamma(p)=p$. For
any $\gamma_i\in \Bbb Z_{q_i}$ such that $\gamma_i\to \gamma$, $o(\gamma_i)=q_i\to \infty$.
\end{remark}

\subsection*{b. Negative curvature and discrete limit isometry groups}
~

A geometric property of a complete metric of negative Ricci curvature is that
if $M$ is compact, then the isometry group is discrete and thus finite (\cite{Bo}). The
discreteness does not hold if $M$ is not compact, e.g., $\dim(\op{Isom} (\Bbb H^n))=\frac{n(n+1)}2$.

In the proof of Theorem B and Theorem D, we need the following property.

\begin{theorem}
Assume an equivariant convergent sequence
satisfying the following commutative diagram:
$$\begin{CD} (\tilde M_i,\tilde p_i,\Gamma_i)@>GH>>(\tilde X,\tilde p,G)
\\ @VV \pi_i V   @VV \pi  V \\
(M_i,p_i)@>GH>> (X,p),\end{CD}$$
where $M_i$ is a compact $n$-manifold of
$\op{diam}(M_i)\le d$, $\Gamma_i=\pi_1(M_i,p_i)$. If $\tilde X$ is isometric to a hyperbolic
manifold, then the identity component $G_0$ is either trivial or not nilpotent.
\end{theorem}

Let $\phi\in \op{Isom}(\Bbb H^n)$. Then $\phi$ acts on the boundary at
infinity of $\Bbb H^n$. From the Poincar\'e model, by Brouwer fixed point theorem
one sees that $\phi$ has a fixed point on the union of $\Bbb H^n$ with its
boundary at infinity. Moreover, $\phi$ satisfies one and only one of the
following property: $\phi$ has a fixed point
in $\Bbb H^n$, $\phi$ has no fixed point in $\Bbb H^n$ and a unique
fixed point or two fixed points on the boundary at infinity; and $\phi$ is
called elliptic, parabolic and hyperbolic respectively (cf. \cite{Ra}).

\begin{lemma}
Let $M$ be a complete non-compact
hyperbolic manifold. Assume that $G$ is a closed group of isometries,
$G_0$ is nilpotent and $M/G$ is compact. Then

\noindent (2.6.1) $G_0$ contains no nontrivial compact subgroup.

\noindent (2.6.2) If $M=\Bbb H^n$, then the center of $G_0$ contains no hyperbolic
element.
\end{lemma}

Note that in Lemma 2.6, $G_0$ may not be trivial; e.g., in the half-plane
model for $\Bbb H^n$, $\op{Isom}(\Bbb H^n)$ contains $\Bbb R^{n-1}$
consisting of parabolic elements which fix the same point $p_\infty$ in the boundary at
infinity. Let $Z=\left<\Bbb R^{n-1},\gamma\right>$, where $\gamma$
is some hyperbolic element which fixes $p_\infty$. Then
$\Bbb H^n/Z$ is a circle. Hence,
to prove Theorem 2.5 i.e., to
rule out parabolic elements in $G_0$, we have to use the fact that
$G$ is the limiting group of an equivariant convergent sequence.

\begin{proof}[Proof of Lemma 2.6]
~

(2.6.1) Since $G_0$ is nilpotent, $G_0$ has a unique
maximal compact subgroup $T^s$ which is also contained in the center $Z(G_0)$ (Lemma 3, \cite{Wi}).
The uniqueness implies that $T^s$ is normal in $G$.  We shall show that $s=0$.

If $s\ge 1$, let $v_1,\dots,v_s$ denote a basis for the lattice $\Bbb Z^s$ ($T^s=\Bbb R^s/\Bbb Z^s$).
Then $H_i=\exp_etv_i$ is a circle subgroup and $T^s=\prod_{i=1}^sH_i$.
The isometric $H_i$-action defines a Killing field $X_i$ on $M$:
$$X_i(x)=\left.\frac{d(H_i(t)(x))}{dt}\right|_{t=0},\quad x\in M.$$
We define a function on $M$ (cf. \cite{Ro1}),
$$f(x)=\frac 12\det(g(X_i,X_j))(x),\quad x\in M.$$
Note that $f(x)$ can be viewed as $\frac 12$-square of the $s$-dimensional volume of $T^s(x)$,
in particular $f(x)$ is independent of the choice of $v_1,\dots,v_s$.

Since $T^s$ is normal in $G$, for $\alpha\in G$, $\alpha(T^s(x))=T^s(\alpha(x))$
and thus $f(\alpha(x))=f(x)$. Since $f$ is $G$-invariant and $M/G$ is compact,
we may assume that $f(x)$ achieves a maximum at $y\in M$, and thus $\Delta f(y)\le 0$.
We claim that $f(x)$ satisfies $\Delta f(y)>0$ at any $y$ such that $f(y)>0$, and
thus a contradiction.

To verify the claim, we first assume that $g_{ij}(y)=g(X_i, X_j)(y)=\delta_{ij}$. Taking any
vector fields $V_1,\dots,V_{n-s}$ on a slice of $T^s(y)$ at $y$ such that $g(V_i,V_j)(y)=\delta_{ij}$
and $g(X_i, V_j)(y)=0$, via the $T^s$-action we extend $V_1, \dots, V_{n-s}$ to be vector fields
on the tube of $T^s(y)$. By construction, $X_1(y),\dots,X_s(y), V_1(y)$,
\dots,$V_{n-s}(y)$ is an orthonormal basis for $T_yM$. For any vector field, $Y$, by calculation we get
$$\begin{aligned} Y(f)(y)&=\frac Y2\left(g_{11}\cdots g_{ss}-\sum_{1\leq i<j\leq s}g_{ij}^2
g_{11}\cdots \hat{g_{ii}}\cdots \hat{g_{jj}}\cdots g_{ss}+ R\right)(y)\\
&=\frac{1}{2}\sum_{i=1}^s g_{11}\cdots
Y(g_{ii})\cdots g_{ss}(y)=\frac{1}{2}\sum_{i=1}^sY(g_{ii})(y),\end{aligned}$$
$$Y(Y(f))(y)=\frac{1}{2}\sum_{i=1}^s Y(Y(g_{ii}))(y)+\sum_{1\leq i<j\leq s}
\left[Y(g_{ii})Y(g_{jj})-(Y(g_{ij}))^2\right](y).$$
Since $[X_i,X_j]=0$ and $X_k(g_{ij})=0$, by calculation we get
$$\begin{aligned} \Delta f(y) = &\sum_{j=1}^s \op{Hess}f(X_j, X_j)(y)+\sum_{l=1}^{n-s}\op{Hess}f(V_l, V_l)(y)
\\
=& \frac{1}{2}\sum_{i=1}^s \Delta g_{ii}(y) +\sum_{l=1}^{n-s}\sum_{1\leq i<j\leq s}
\left[V_l(g_{ii})V_l(g_{jj})-(V_l(g_{ij}))^2\right](y).
\end{aligned}$$
Since for any vector fields $V, W$, any $1\leq k\leq s$, $g(\nabla_V X_k, W)= -
g(\nabla_W X_k, V)$,
$$\begin{cases} \frac{1}{2}\Delta g_{ii}(y)= \sum_{j=1}^s|\nabla_{X_j}X_i|^2(y)+
\sum_{l=1}^{n-s}|\nabla_{V_l}X_i|^2(y)-\op{Ric}(X_i, X_i)(y)\\
|\nabla_{X_j}X_i|^2(y)=\sum_{k=1}^s g^2(\nabla_{X_j}X_i, X_k)(y)+\sum_{l=1}^{n-s} g^2(\nabla_{X_j}X_i, V_l)(y)\\
|\nabla_{V_l}X_i|^2(y)=\sum_{k=1}^s g^2(\nabla_{V_l}X_i, X_k)(y)+\sum_{k=1}^{n-s} g^2(\nabla_{V_l}X_i, V_k)(y).\end{cases}$$
Finally,
$$\begin{aligned} \Delta f(y)=& 2\sum_{l=1}^{n-s}\left[\sum_{i=1}^s g(\nabla_{V_l}X_i, X_i)(y)\right]^2+ \sum_{i,j,k=1}^s g^2(\nabla_{X_j}X_i, X_k)(y)\\
&+\sum_{i=1}^s\sum_{k,l=1}^{n-s}g^2(\nabla_{V_l}X_i, V_k)(y)-\sum_{i=1}^s \op{Ric}(X_i, X_i)(y).\end{aligned}$$
In particular, we conclude that if $f(y)>0$ i.e., $X_1(y),...,X_s(y)$ are linear independent, then $\Delta f(y)>0$.

In general, at $y$ where $f(y)>0$ we may choose Killing vector fields, $W_1,...W_s$, such that $W_1(y),...,W_s(y)$ is orthonormal at $y$, and
let $A=(a_{ij})$ be a constant $n\times n$-matrix such that $W_i(y)=\sum_{j=1}^sa_{ij}X_j(y)$. Then $f(x)=\frac 12\det(AA^T)\cdot \det(g(W_i,W_j))(x)$,
and thus $\Delta f(y)>0$ at $y$ where $f(y)>0$.

(2.6.2) Since $G_0$ is nilpotent, by (2.6.1) we may assume that $Z(G_0)=\Bbb R^s$ is not
trivial i.e., $s\ge 1$. Assume that $\phi\in Z(G_0)$ is a hyperbolic element i.e., $\phi$ acts
freely on $\Bbb H^n$ and has two fixed points on the boundary at infinity. Let
$c(t)$ be the unique minimal geodesic connecting the two $\phi$-fixed points.
Then $\phi$ preserves $c(t)$, and $c(t)$ is the unique line in $\Bbb H^n$
preserved by $\phi$ (because if a line $\alpha(t)$ is preserved by $\phi$, then
$c(t)$ and $\alpha(t)$ are preserved by $\phi^2$ which fixes the
two ends). Since any element in $G_0$ commutes with
$\phi$, $G_0$ preserves $c(t)$, and thus $G_0=Z(G_0)=\Bbb R^1$ such that
$c(t)$ is an $\Bbb R^1$-orbit, which is the unique line $\Bbb R^1$-orbit.
Since $\Bbb R^1$ is normal in $G$, any element in $G$ preserves $c(t)$, and thus
$G/\Bbb R^1$ has a fixed point on $\Bbb H^n/\Bbb R^1$. Since $G/\Bbb R^1$ is discrete,
$G/\Bbb R^1$ is finite. On the other hand, $\Bbb H^n/\Bbb R^1$ is not
compact, because otherwise for $\Bbb Z\subset \Bbb R^1$, $\Bbb H^n/\Bbb Z$
is compact hyperbolic manifold on which $\Bbb R/\Bbb Z$ acts isometrically,
a contradiction. Since $\Bbb H^n/\Bbb R^1$ is not compact and $G/\Bbb R^1$ is finite,
$\Bbb H^n/G=(\Bbb H^n/\Bbb R^1)/(G/\Bbb R^1)$ is not compact, a contradiction.
\end{proof}

\begin{proof}[Proof of Theorem 2.5]
~

Assume that $G_0$ is nilpotent. We shall show that $G_0=e$.

By (2.6.1), we assume that $G_0$ acts freely on $\tilde X$.  We first assume that
$\tilde X=\Bbb H^n$. By (2.6.2), $G_0$ contains only parabolic elements.
Since $G_0$ is parabolic, in the upper half plane model we see that $G_0(\tilde p)$
is contained in the horizontal hyperplane $\Bbb R^{n-1}$. Since $\Bbb R^{n-1}$
contains no segment, any $G_0$-orbit contains no piece of minimal geodesic.
We shall derive a contradiction by constructing a sequence of minimal geodesic
$\gamma_i$ on $\tilde M_i$ that converges to a minimal geodesic in some $G_0$-orbit.

Let $v\in T_eG_0$ be a unit vector, let $\phi=\exp_ev$. Let $t_k=\frac 1k\in [0,1]$,
and let $\phi_k=\exp_et_kv\in G_0$. From the equivariant convergent commutative
diagram,
$$\begin{CD} (\tilde M_i,\tilde p_i,\Gamma_i)@> GH >>(\Bbb H^n,\tilde p,G)
\\@VV \pi_i V   @VV \pi V \\
(M_i,p_i)@> GH >>(X,p),
\end{CD}$$
we may assume $\gamma_{i,k}\in \Gamma_i$ such that $\gamma_{i,k}\to \phi_k$, and thus
for any $1\le j\le k$, $\gamma_{i,k}^j\to \phi_k^j$. Since $M_i$ is compact, we may
assume that $p_{i,k}$ is chosen so that $\gamma_{i,k}$ is represented by a close
geodesic $c_{i,k}$ at $p_{i,k}$. Consequently, the lifting $\tilde c^k_{i,k}$ of $c^k_{i,k}(t)$ at $\tilde p_{i,k}$
is a segment that contains a piece of length almost one. Let $\tilde c^k_{i,k}\to \tilde c_k\subset \Bbb H^n$.
Clearly, $\tilde c_k$ is a segment. Let $k\to \infty$ and via a standard diagonal argument
we conclude that $\tilde c_k\to \tilde c$ is contained in $G_0(\tilde p)$.

If $\tilde X\ne \Bbb H^n$, we consider the lifting isometric
$G_0$-action on $\Bbb H^n$ satisfying
the following diagram commutes:
$$\begin{CD} G_0\times \Bbb H^n@> \tilde \mu>>\Bbb H^n
\\@VV \op{id}\times \pi V   @VV \pi V \\
G_0\times \tilde X@> \mu >>\tilde X.
\end{CD}$$

If $Z(G_0)$ contains a parabolic element, then following the above argument
we see that $G_0$-orbit in $\tilde X$ contains a piece of minimal
geodesic, and thus its lifting to $\Bbb H^n$ is a piece of minimal geodesic in a $G_0$-orbit in
$\Bbb H^n$, a contradiction.

If $Z(G_0)$ contains a hyperbolic element, then by the proof of (2.6.2) we see that
$G_0=\Bbb R^1$ and $G/G_0$ fixes a point in $\tilde X/\Bbb R^1$ (note that $\pi_1(\tilde X)$
commutes with the lifting $G_0$-action), which contradicts
to that $\tilde X/G$ is compact.
\end{proof}

\subsection*{c. The $C^0$-convergence}
~

In the proof of Theorem A, the following $C^0$-convergence plays an important role
(see the proof of (3.5.2)). Let $(M,g)$ be a compact Riemannian manifold, and let
$g(t)$ denote the Ricci flow i.e., the solution of the following PDE (\cite{Ha1}):
$$\frac{\partial g(t)}{\partial t}=-2\op{Ric}(g(t)),\qquad g(0)=g.$$

\begin{theorem}
Let $g_i$ ($i=0, 1$) be two Riemannian metrics
on a compact $n$-manifold $M$ such that $\op{Ric}_{g_1}\ge -(n-1)$. Given
$\epsilon>0$, there are constants, $\delta(\epsilon,g_0), T=T(n,\epsilon,g_0)>0$,
such that for $0<\delta\le \delta(\epsilon,g_0)$, if
$$\op{id}_M:(M, g_1)\to (M,g_0)\text{ is a $\delta$-GHA},$$
then the Ricci flow $g_1(t)$ exists for all $t\in (0,T]$
such that $|g_1(T)-g_0|_{C^0(M)}<\epsilon$.	
\end{theorem}

Note that the existence of $T(n,\epsilon,g_0)$ is a consequence of
the Perel'man's pseudo-locality (Theorem 10.1, Corollary 10.2 in \cite{Pe2}).
For our purpose, we state it in the following form (\cite{CM}, \cite{TW}).

\begin{theorem}
Given $n, \delta>0$, there exist constants, $r(n),
\epsilon(n,\delta), C(n)$, $T(n,\delta)>0$, such that if a compact $n$-manifold $(M,g)$ satisfies
$$\op{Ric}_g\ge- (n-1),\quad d_{GH}(B_r(x),\sball{0}{r})<\epsilon(n,\delta) r,\quad 0<r<r(n), \,\, x\in M,$$
then the Ricci flow $g(t)$ exists for all $t\in [0,T(n,\delta)]$ and satisfies
$$|\op{Rm}(g(t))|_{M}\le \frac \delta t,\quad \op{vol}(B_{\sqrt t}(x,g(t)))
\ge C(n)(\sqrt t)^n.$$
\end{theorem}

By (1.7.2), a sequence of compact $n$-manifolds, $M_i\overset{GH}{\longrightarrow}M$, such that $\op{Ric}_{M_i}\ge -(n-1)$
and $M$ is a Riemannian $n$-manifold is equivalent to a sequence of Riemannian metrics on $M$,
$g_i$ and $g$, such that $\op{id}_M: (M,g_i)\to (M,g)$ is an $\epsilon_i$-GHA, $\epsilon_i\to 0$.

\begin{corollary}
Assume a sequence of Riemannian metrics, $g_i$,
and a Riemannian metric $g$ on a compact $n$-manifold $M$ satisfying
$$\op{Ric}_{g_i}\ge -(n-1),\quad \op{id}_M:(M, g_i)\to (M,g)
\text{ is an $\epsilon_i$-GHA}, \quad \epsilon_i\to 0.$$
Then passing to a subsequence there is a sequence of Ricci
flow solutions of $g_i$ at time $t_i\to 0$, $g_i(t_i)$,
such that $|g_i(t_i)-g|_{C^0(M)}\to 0$
as $i\to \infty$.
\end{corollary}

In the proof of Theorem 2.7, we need the following property for the distance
function of $g(t)$, which is due to Bamler-Wilking (\cite{BW}).

\begin{lemma}
Let the assumption be as in Theorem 2.8. There exists $0<\eta(n,\delta)<T(n,\delta)$ such that for any $x, y\in M$
with $|xy|_{g(t)}<\sqrt t\le \eta(n,\delta)$,
$$||xy|_g-|xy|_{g(t)}|\le \Psi(\delta|n)\sqrt t.$$
\end{lemma}

\begin{proof}
	Because $g(t)$ satisfies that $\op{Ric}_{g(t)}\le \frac {(n-1)\delta} t$,
	it is known that the function, $|xy|_{g(t)}+25(n-1)\sqrt{\delta t}$, is monotonically
	increasing in $t$ (cf. 17. of \cite{Ha2}, Corollary 3.26 in \cite{MT}).
	Consequently, $|xy|_{g(t)}+25(n-1)\sqrt{\delta t}\ge |xy|_g$.
	
	To prove an opposite inequality, we will assume that $|xy|_{g(t)}<\sqrt t$. By Theorem 2.8 and
	the injectivity radius estimate, we may assume that $\op{injrad}(x,g(t))\ge \rho\sqrt t$ for
	all $x$, where $\rho$ is a constant depending on $n$. Without loss of generality we may assume that $\rho\ge 1$.
	
	Arguing by contradiction, assume some $\sigma>0$ and given any $\delta_i\to 0$, there is a sequence
	of compact $n$-manifolds $(M_i,g_i)$, $x_i, y_i\in M_i$ and $t_i\in (0,T(n,\delta_i)]$ with $t_i\to 0$ such that
	$|x_iy_i|_{g_i(t_i)}>|x_iy_i|_{g_i}+\sigma\sqrt {t_i}$.
	Let $d_i=|x_iy_i|_{g_i(t_i)}$. It is easy to check the following relations
	(assume that $25(n-1)\sqrt{\delta_i}<\frac \sigma 4$):
	$$\begin{cases} B_{d_i-25(n-1)\sqrt{\delta_i t_i}-\frac \sigma 2\sqrt{t_i}}(x_i,g_i(t_i))
	\subset B_{d_i-\frac \sigma2\sqrt{t_i}}(x_i,g_i)
	\\ B_{\frac \sigma4\sqrt{t_i}}(y_i,g_i(t_i))\subset
	B_{d_i-\frac \sigma2\sqrt{t_i}}(x_i,g_i)\end{cases}$$
	Let $\ell_i=\frac {d_i}{\sqrt {t_i}}$, and let $s_i=25(n-1)\sqrt{\delta_i}-\frac \sigma 2$. Then $\sigma<\ell_i\le 1$
	and $s_i\to -\frac \sigma2$.
	Since $B_{\frac \sigma4\sqrt{t_i}}(y_i,g_i(t_i))\cap B_{d_i-s_i\sqrt{t_i}}(x_i,g_i(t_i))=\emptyset$,
	by (\cite{Ha1}) and Bishop-Gromov volume comparison we derive
	$$\begin{aligned} &\frac{\svolball{-t_i}{\ell_i-\frac \sigma2}}{(\sqrt{t_i})^n}
	=\svolball{-1}{d_i-\frac\sigma2\sqrt{t_i}}\ge \op{vol}_{g_i}(B_{d_i-\frac\sigma2\sqrt{t_i}}(x_i,g_i))\\
	&\ge \op{vol}_{g_i}(B_{d_i-s_i
		\sqrt{t_i}}(x_i,g_i(t_i)))+\op{vol}_{g_i}(B_{\frac \sigma4\sqrt{t_i}}
	(y_i,g_i(t_i)))\\
	&\ge (1-\Psi(t_i|n))\left[\op{vol}_{g_i(t_i)}(B_{d_i-s_i\sqrt{t_i}}(x_i,g_i(t_i)))+\op{vol}
	_{g_i(t_i)}(B_{\frac \sigma4\sqrt{t_i}}(y_i,g_i(t_i)))\right]
	\\
	&\ge (1-\Psi(t_i|n))\left(\frac {\svolball
		{\delta_i}{\ell_i-s_i}}
	{(\sqrt{t_i})^n}+\frac{\svolball{\delta_i}{\frac
			\sigma4}}{(\sqrt{t_i})^n}\right),
	\end{aligned}$$
	where the last inequality is because $\text{sec}_{t_i^{-1}g_i(t)}\le \delta_i$ and
	$\op{injrad}(x_i,g_i(t))\ge \rho\sqrt t$.
	We may assume that $\ell_i\to \ell$, $\sigma\le \ell\le 1$. As $i\to \infty$,
	from the above we conclude that $\svolball{0}{\ell-\frac \sigma2}\ge
	\svolball{0}{\ell-\frac \sigma2}+\svolball{0}{\frac \sigma4}$, a contradiction.
\end{proof}

\begin{remark} 
    Note that $\eta(n,\delta)\to 0$ as $\delta\to 0$, and it is unlikely that $T(n,\delta)\to 0$ as $\delta\to 0$.
\end{remark}

\begin{proof}[Proof of Theorem 2.7]
	~
	
	Let $\op{id}_M: (M,d_{g_1})\to (M,d_{g_0})$ be a $\delta$-GHA, where $\delta$
	will be specified later. By Theorem 1.6, given $\delta_1>0$, we may assume $\delta$ small
	so that $(M,g_1)$ satisfies the conditions of Theorem 2.8 with $\epsilon(n,\delta_1)$ and $r=r(g_0)$,
	and thus there are constants, $C(n), T=T(n,\delta_1,g_0)>0$, such that
	the Ricci flow solution $g_1(t)$ with $t\in (0,T]$ satisfies that
	$$|\op{Rm}(g_1(t))|_{M}\le \frac {\delta_1} t,\quad \op{vol}(B_{\sqrt t}(x,g_1(t)))\ge C(n)(\sqrt t)^n.$$
	For all $x\in M$, the re-scaling metric satisfies that
	$$|\op{Rm}(T^{-1}g_1(T))|_{M}\le \delta_1,\quad \op{vol}(B_1(x,T^{-1}g_1(T)))\ge C(n).$$
	By Lemma 2.10,
	$$\op{id}_{B_1(x,T^{-1}g_1)}: (B_1(x,T^{-1}g_1),d_{T^{-1}g_1})\to (B_1(x,T^{-1}g_1),d_{T^{-1}g_1(T)})$$
	is an $\Psi(\delta_1|n)$-GHA, and thus
	$$\op{id}_{B_{\frac 12}(x,T^{-1}g_0)}: (B_{\frac 12}(x,T^{-1}g_0),d_{T^{-1}g_1(T)})\to (B_{\frac 12}(x,T^{-1}g_0),d_{T^{-1}g_0})$$
	is an $(\Psi(\delta_1|n)+\frac \delta T)$-GHA. By Cheeger-Gromov $C^{1,\alpha}$-convergent theorem (cf. \cite{Pet}),
	we first choose $\delta_1=\delta_1(\epsilon,g_0)$ small so that $\op{id}_{B_1(x,T^{-1}g_0)}$ is an $2\Psi(\delta_1|n)$-GHA
	implies that $|T^{-1}g_1(T)-T^{-1}g_0|_{C^{1,\alpha}(B_{\frac 12}(x,T^{-1}g_0))}<\epsilon$.
	Note that $\op{id}_{B_1(x,T^{-1}g_0)}$ is an $2\Psi(\delta_1|n)$-GHA if we choose $\delta(\epsilon,g_0)=\Psi(\delta_1|n)\cdot T$.
	Since the $C^0$-norm is scaling invariant, $|g_1(T)-g_0|_{C^0(B_{\frac {\sqrt T}2}(x,g_0))}$ $<\epsilon$.
	Because $x\in M$ is arbitrary, $|g_1(T)-g_0|_{C^0(M)}<\epsilon$.
\end{proof}

\section{Proofs of Theorem A-C, Theorem E and Theorem 0.4}

Consider a sequence of compact $n$-manifolds, $M_i\overset{GH}{\longrightarrow}X$, $\epsilon_i\to 0$,
$$\op{Ric}_{M_i}\ge (n-1)H,\, \op{diam}(M_i)\le d,\,\volume(B_1(\tilde p_i))\ge v,\,\frac{\volume(B_\rho(x_i^*))}
{\svolball H\rho}\ge 1-\epsilon_i. \eqno (3.1.1)$$
From Section 1, subsection b, passing to a subsequence we may assume
the following commutative diagram:
$$\begin{CD} (\tilde M_i,\tilde p_i,\Gamma_i)@>GH>>(\tilde X,\tilde p,G)
\\@VV \pi_i V   @VV \pi  V \\
(M_i,p_i)@>GH>>(X,p),
\end{CD}\eqno (3.1.2)$$
where $\Gamma_i$ denotes the deck transformation group, $G$ is a closed subgroup of $\op{Isom}(\tilde X)$,
which is a Lie group (Theorem 1.8).

\addtocounter{theorem}{1} 
\begin{lemma}
	Let $\tilde X$ be as in the above. Then
	$\tilde X$ is isometric to Riemannian $n$-manifold of constant curvature $H$.
\end{lemma}

\begin{proof}
	For $\tilde x\in \tilde X$, let $\tilde x_i\in \tilde M_i$
	such that $\tilde x_i\to \tilde x$.
	Let $x_i=\pi_i(\tilde x_i)$,  and let $\pi_i^*: (U^*_\rho(x_i^*),x_i^*)
	\to (B_\rho(x_i),x_i)$ be the Riemannian universal covering.
	Consider the commutative diagram in (0.17), and by Theorem 1.2,
	$$d_{GH}(B_{\frac \rho2}(x^*),\sball H{\frac \rho2})=
	\lim_{i\to \infty}d_{GH}(B_{\frac \rho2}(x_i^*),\sball H{\frac \rho2})\le \lim_{i\to
		\infty}\Psi(\epsilon_i|n,\rho,H)=0,$$
	and thus $B_{\frac \rho2}(x^*)$ is isometric to $\sball H{\frac \rho2}$. By Bishop-Gromov
	relative volume comparison, the condition $\op{vol}(B_1(\tilde p_i))\ge v$ implies that
	for any $\tilde x_i\in \tilde M_i$, $\op{vol}(B_\rho(\tilde x_i))
	\ge v(n,\rho,d,H, v)>0$. By Corollary 2.2, we can conclude that
	$K$ acts freely on $B_{\frac \rho4}(x^*)$, and thus $B_{\frac \rho4}(\tilde x)$ is
	a manifold of constant curvature $H$. Consequently, $\tilde X$ is a manifold of constant curvature $H$.
\end{proof}

\begin{corollary}
	Let the assumptions be as in Theorems A-C (resp.
	$H=1, -1$ or $0$).
	Then there is $\rho'(n,\rho,d,v)>0$ such that the Riemannian
	universal covering $\tilde M$ satisfies
	$$\frac{\volume(B_{\rho'}(\tilde x))}{\svolball{H}{\rho'}}\ge 1-\Psi
	(\epsilon|n,\rho,d,v),\qquad \tilde x\in \tilde M.$$
\end{corollary}

\begin{proof}
	Arguing by contradiction, assume $\rho_k\to 0$ such
	that for each $\rho_k$ there is $\epsilon(\rho_k)>0$ and a sequence
	$M_{i,k}$ such that
	$$\frac{\volume(B_\rho(x_{i,k}^*))}{\svolball H\rho}\ge 1-
	\epsilon_i\to 1, \quad \forall\, x_{i,k}\in M_{i,k},$$
	and there is $\tilde q_{i,k}\in \tilde M_{i,k}$ such that
	$$\frac{\volume(B_{\rho_k}(\tilde q_{i,k}))}{\svolball{H}
		{\rho_k}}< 1-\epsilon(\rho_k), \quad \forall \, i.
	\eqno (3.3.1)$$
	Passing to a subsequence, we may assume
	$$\begin{CD}(\tilde M_{i,k},\tilde q_{i,k})@>GH>>(\tilde X_k,\tilde q_k)\end{CD}.$$
	By Lemma 3.2, $\tilde X_k$ is isometric to space form of constant curvature
	$H$ and $\op{vol}(B_1(\tilde q_k))\ge v'(n,d,v)>0$ (Theorem 1.6).
	By Cheeger's injectivity estimate, we may assume that
	$\op{injrad}(\tilde q_k)\ge \rho'(n,\rho,d,v)>0$.
	For fixed $\rho_k<\frac{\rho'}2$, by Theorem 1.6 we have that
	$\op{vol}(B_{\rho_k}(\tilde q_{i,k}))\to \svolball{H}{\rho_k}$,
	a contradiction to (3.3.1).
\end{proof}

\subsection*{a. Proofs of Theorem A-C}
~

Consider a sequence in (3.1.1) and (3.1.2) with $H=1$, and thus
$\tilde X$ is isometric to $S^n_1$ (Lemma 3.2, Theorem 1.7). In
the proof of Theorem A, we need the following result in \cite{MRW}.

\begin{lemma}
	Let $M_i\overset{GH}{\longrightarrow}X$ be a sequence of compact
	$n$-manifolds satisfying
	$$\op{Ric}_{M_i}\ge -(n-1),\quad \op{diam}(M_i)\le d,\quad \op{vol}
	(B_1(\tilde p_i))\ge v>0,$$
	and the commutative diagram (3.1.2). If $\Gamma_i$ is finite,
	then for $i$ large, there is an
	injective homomorphism, $\phi_i: \Gamma_i\to G$, which is also
	an $\epsilon_i$-GHA with $\epsilon_i\to 0$.
\end{lemma}

Note that Lemma 3.4 was originally stated in \cite{MRW} under the condition that $\op{sec}_{M_i}\ge -1$. Because
the sectional curvature condition was used only to conclude that a limiting group is a Lie group,
by Theorem 1.8 Lemma 3.4 is valid when `$\op{sec}_{M_i}\ge -1$' is replaced by `$\op{Ric}_{M_i}\ge -(n-1)$'.

Let $\phi_i: \Gamma_i\to G$ be as in Lemma 3.4.  By Theorem 1.7, we may
assume a diffeomorphism, $\tilde h_i: \tilde M_i\to S^n_1$, such that
$(\tilde h_i,\phi_i)$ is also an $\epsilon_i$-equivariant
GHA i.e., for all $\tilde x_i\in \tilde M_i$ and $\gamma_i\in \Gamma_i$,
$$|\tilde h_i(\tilde x_i)[\phi_i(\gamma_i)\tilde h_i(\gamma_i^{-1}(\tilde x_i))]|<\epsilon_i.$$
Note that via $\tilde h_i$, $\Gamma_i$ acts freely on $\tilde X$:
$\gamma_i(\tilde x)=\tilde h_i(\gamma_i(\tilde h_i^{-1}(\tilde x)))$ for $\tilde x\in \tilde X$ and
$\gamma_i\in \Gamma_i$.  We shall use $\Gamma_i(\tilde h_i)$ to denote the $\Gamma_i$-action on $\tilde X$
via $\tilde h_i$.

\begin{theorem}
	Let $M_i$ be a sequence of compact $n$-manifolds satisfying
	$$\op{Ric}_{M_i}\ge (n-1),\quad \frac{\op{vol}_\rho(B_\rho(\tilde x_i))}{\svolball{1}{\rho}}
	\ge 1-\epsilon_i\to 1,\quad \tilde x_i\in \tilde M_i$$
	and the commutative diagram (3.1.2). Then for $i$ large,
	
	\noindent (3.5.1) $\phi_i(\Gamma_i)$ acts freely on $S^n_1$.
	
	\noindent (3.5.2) The $\Gamma_i(\tilde h_i)$-action and the $\phi_i(\Gamma_i)$-action on $S^n_1$ are conjugate.
\end{theorem}

\begin{proof}
(3.5.1)  If $e\ne \gamma_i\in \Gamma_i$, $\tilde y\in S^n_1$ such that $\phi_i(\gamma_i)(\tilde y)=\tilde y$,
then $\left<\gamma_i\right>\to \Lambda\ne e$ (Lemma 2.3) and $\Lambda(\tilde y)=\tilde y$.  Without loss
of generality, we may assume $\tilde y$ is chosen such that $\tilde x_i\to \tilde y$ and
the displacement of $\gamma_i$ achieves a minimum at $\tilde x_i$.  Since $\left<\gamma_i\right>(\tilde x_i)\overset{GH}{\longrightarrow}
\Lambda (\tilde y)=\tilde y$,
$r_i=\op{diam}(\left<\gamma_i\right>(\tilde x_i))\to 0$. Consider the rescaling sequence,
$$\begin{CD}(r_i^{-1}\tilde M_i,\tilde x_i,\left<\gamma_i\right>)@>GH>>(\Bbb R^n,v,K).\end{CD}$$
Since $\op{diam}K(v)=1$, $K$ is compact.  Then $K$ has a
fixed point, say $0$, and let $\tilde z_i\in r_i^{-1}\tilde M_i$ such that $\tilde z_i\to 0$.
Then $\left<\gamma_i\right>(\tilde z_i)\to K(0)=0$. This is not possible, because
$$\op{diam}(\left<\gamma_i\right>(\tilde z_i))\ge \op{diam}(\left<\gamma_i\right>(\tilde x_i))=1,$$
a contradiction.

(3.5.2) Let $\tilde g_i$ denote the pullback metric on $S^n$ by $\tilde h_i^{-1}$. Then
the identity map, $\op{id}_{S^n}: (S^n,\tilde g_i,\Gamma_i(\tilde h_i))\to (S^n, \b g^1,\phi_i(\Gamma_i))$,
is an $\epsilon_i$-equivariant GHA. Following \cite{GK}, we will construct an equivariant map
via the method of center of mass with respect to $\b g^1$: fixing $\tilde x\in S^n$, let
$A(\tilde x)=\{\phi_i(\gamma_i)^{-1}(\gamma_i(\tilde x)),\,\,\gamma_i\in \Gamma_i(\tilde h_i)\}$.
Since $A(\tilde x)\subset B_{\frac \pi4}(\tilde x)$, $A(\tilde x)$ has a center of mass,
say $\tilde y$.  We then define $\tilde  f_i: S^n_1\to S^n_1$
by $\tilde f_i(\tilde x)=\tilde  y$. Then $\tilde f_i$ is a differentiable map satisfying
that $\tilde f_i(\gamma_i(\tilde x))=\phi_i(\gamma_i)(\tilde f_i(\tilde x))$.

According to \cite{GK}, $\tilde f_i$ is a diffeomorphism if the two actions are
$\epsilon$-close in $C^1$-norm i.e.,
$$\max\{|\tilde x\phi_i(\gamma_i)^{-1}\gamma_i
(\tilde x)|_{\b g^1},\, \tilde x\in S^n\}<\epsilon, \,
|d(\phi_i(\gamma_i)^{-1}\gamma_i)(X)-\b P(X)|_{\b g^1}<\Psi(\epsilon),$$
for all $\gamma_i
\in \Gamma_i(\tilde h_i)$ and $|X|_{\b g^1}=1$, where $\b P$ denotes the
$\b g^1$-parallel translation along the unique minimal geodesic joining $\tilde x$ and
$\phi_i(\gamma_i)^{-1}\gamma_i(\tilde x)$ and $\epsilon>0$ is a constant determined by $\b g^1$.

Given $\epsilon>0$, by Theorem 2.7 we may assume that $\op{id}_{S^n}: (S^n,\tilde g_i(T))\to (S^n,\b g^1)$ is an
$\delta(\epsilon)$-GHA for $i$ large, where $T=T(n,\epsilon,\b g^1)>0$
such that $|T^{-1}\tilde g_i(T)-T^{-1}\b g^1|_{C^{1,\alpha}
	(B_{\frac 12}(\tilde x,T^{-1}\b g^1))}<\epsilon$ (see the end of proof of Theorem 2.7).
Consequently, restricting to $B_{\frac
	12}(\tilde x,T^{-1}\b g^1)$, exponential maps of $T^{-1}g_i(T)$ and $T^{-1}\b g^1$ are
$C^\alpha$-close, and therefore the $\Gamma_i(\tilde h_i)$ and $\phi_i(\Gamma_i)$-actions are
$\epsilon$-close in $C^1$-norm. Since $\epsilon>0$ is arbitrary, the desired conclusion
follows.
\end{proof}

\begin{proof}[Proof of Theorem A]
	~
	
	Arguing by contradiction, assume a sequence, $M_i\overset{GH}{\longrightarrow}X$, satisfying (3.1.1) and (3.1.2)
	for $H=1$ such that $M_i$ is not diffeomorphic to any spherical $n$-space form.  By Lemma 3.2,
	$\tilde X$ is isometric to spherical space form. By Theorem 1.7, $\tilde X$ is diffeomorphic
	to $\tilde M_i$ which is simply connected, and therefore $\tilde X=S^n_1$. By (3.5.1) and
	(3.5.2), $M_i=\tilde M_i/\Gamma_i$ is diffeomorphic to $S^n_1/\phi_i(\Gamma_i)$, a contradiction.
\end{proof}

\begin{proof}[Proof of Theorem B]
	~
	
	Arguing by contradiction, assume a sequence, $M_i\overset{GH}{\longrightarrow}X$, satisfying (3.1.1) and (3.1.2)
	for $H=-1$ such that $M_i$ is not diffeomorphic to any hyperbolic $n$-manifold.  By Lemma 3.2,
	$\tilde X$ is isometric to a hyperbolic $n$-manifold (we do not yet know that $\tilde X$ is simply
	connected).  We claim that there is a constant $c(n,\rho,d,v)>0$ such that $\op{vol}(M_i)\ge c(n,\rho,d,v)$.
	Consequently, $G$ is discrete. By Corollary 3.3 we are able to apply Theorem 2.1 and
	conclude that $G$ acts freely on $\tilde X$ and thus $X=\tilde X/G$ is isometric to a
	hyperbolic $n$-manifold. By Theorem 1.7, $M_i$ is diffeomorphic to $\tilde X/G$, a contradiction.
	
	If the above claim fails, then $\dim(X)<n$ and thus $\dim(G_0)>0$. By Lemma 1.13 there
	is $\epsilon>0$ such that $\Gamma_i^\epsilon\to G_0$.
	By Theorem 1.9, $\Gamma_i^\epsilon$ has a nilpotent subgroup of bounded index, and thus
	$G_0\ne e$ is nilpotent, a contradiction to Theorem 2.5.
\end{proof}

\begin{proof}[Proof of Theorem C]
	~
	
	Arguing by contradiction, we may assume a sequence $M_i\overset{GH}{\longrightarrow}X$ satisfying
	(3.1.1) and (3.1.2) for $H=0$ and $M_i$ is not flat.
	By Lemma 3.2, $\tilde X$ is a flat manifold, and thus $\tilde X=\Bbb R^k\times F^{n-k}$
	and $F^{n-k}$ is a compact flat manifold.
	On the other hand, by Splitting theorem of Cheeger-Gromoll, $\tilde M_i=
	\Bbb R^{k_i}\times N_i$, where $N_i$ is a compact simply connected manifold
	of non-negative Ricci curvature.
	
	We claim that $\op{diam}(N_i)\le D(n)$ a constant depending on $n$,  and without loss of generality
	we may further assume that $\op{diam}(F^{n-k})\le D(n)$. Consequently, for any $R>D(n)$ and $i$ large,
	$B_R(\tilde p_i)$ is simply connected and is diffeomorphic to
	$B_R(\tilde p)$ (Theorem 1.7), which implies that $n-k=0$, and
	thus $N_i$ is a point i.e., $M_i$ is a flat manifold, a contradiction.
	
	Assuming that  $\op{diam}(N_i)=r_i\to \infty$, passing to a subsequence we may assume
	$$\begin{CD} (r_i^{-1}\Bbb R^{k_i}\times N_i,\tilde p_i,\Gamma_i)@> GH>> (\Bbb R^k\times N,\tilde p',G')
	\\  @ VV \pi_i V @VV \pi V\\ (r_i^{-1}M_i, p_i) @> GH>> p,\end{CD}$$
	where $N$ is a compact length space of diameter $1$. Note that $G'=G'_0$
	is a nilpotent group (Theorem 1.9) acting effectively  and transitively on $\Bbb R^k\times N$.
	Consequently, $N$ is a $s$-torus ($s\ge 1$). Since $r_i^{-1}N_i\overset{GH}{\longrightarrow} N=T^s$, there is an onto map
	from $\pi_1(N_i)\to \pi_1(T^s)$ (cf. \cite{Tu}),
	a contradiction \footnote{The proof of $\op{diam}(N_i)\le D(n)$ was due to J. Pan.\hfill{$\,$}}.
\end{proof}

\subsection*{b. Proof of Theorem E}
~

\begin{lemma}
	Given $n, \rho>0$, there exists a constant $\epsilon(n,\rho)>0$ such that for any
	$0<\epsilon<\epsilon(n,\rho)$, if a compact Einstein $n$-manifold $M$ of Ricci curvature $\equiv H$
	satisfies
	$$\frac{\op{vol}(B_\rho(x^*))}{\svolball{H}{\rho}}\ge 1-\epsilon,\quad \forall\, x\in M,$$
	then the sectional curvature is almost constant i.e.,
	$$H-\Psi(\epsilon|n,\rho) \le \op{sec}_M\le H+\Psi(\epsilon|n,\rho).$$
\end{lemma}

\begin{proof}
	Arguing by contradiction, assuming a sequence $\epsilon_i\to 0$
	and a sequence of Einstein $n$-manifolds $M_i$ which satisfy the conditions
	of Lemma 3.6 with respect to $\epsilon_i$, but there are $p_i\in M_i$ and a
	plane $\Sigma_i\subset T_{p_i}M_i$ such that $|\op{sec}(\Sigma_i)-H|\ge \delta>0$.
	
	By Theorem 1.2, passing to a subsequence we may assume that
	$B_\rho(p^*_i)\overset{GH}{\longrightarrow}\sball{H}{\rho}$.  Since for $i$ large, $B_{\frac \rho2}(p_i^*)$ is diffeomorphic
	to $\sball{H}{\frac \rho2}$ (compare to  (1.7.2)),  we may identify the sequence as a sequence of metrics $d_i^*$
	on $\sball{H}{\frac \rho2}$ that converges to $\b d^H$. Since the lifting metrics $g_i^*$ on $B_\rho(p_i^*)$ is Einstein,
	passing to a subsequence we may assume that $g_i^*\overset{C^k}{\longrightarrow} \b g^H$ for any $k<\infty$ (\cite{Ch}).
	In particular, $\op{sec}_{g_i^*}|_{B_{\frac\rho2}(p_i^*)}\to H$ i.e.,
	$$H-\Psi(\epsilon_i|n,\rho) \le \op{sec}_{B_{\frac \rho2}(p_i)}\le H+\Psi(\epsilon_i|n,\rho),$$
	a contradiction.
\end{proof}

\begin{proof}[Proof of Theorem E]
	By Lemma 3.6, $M_i$ has almost constant sectional curvature $H$.
	
	Case 1. Assume $H=-1$. Since $M$ has bounded negative sectional
	curvature, by Heintze-Margulis lemma (\cite{He}) we may assume
	$\op{vol}(M)\ge v(n)>0$. By now the desired conclusion follows
	from Theorem B.
	
	Case 2. Assume $H=0$. Then $M$ is almost flat, and thus by Gromov's almost
	flat manifolds theorem $\tilde M$ is contractible (\cite{Gr}). By Cheeger-Gromoll Splitting
	theorem it follows that $M$ is flat.
	
	Case 3. Assume $H=1$. First, since the curvature is almost one, the classical $1/4$-pinched
	injectivity radius estimate implies that $\tilde M$ has injectivity radius $> \frac \pi2$.
	By now the desired conclusion follows from Theorem A.
\end{proof}

\begin{remark}
	In a recent paper \cite{CRX}, we generalized Theorem E to manifolds with bounded Ricci curvature.
\end{remark}

\subsection*{c. Proof of Theorem 0.4}
~

We first extend Theorem C to a limit space.

\begin{lemma}
	Given $n, \rho,v>0$, there is $\epsilon_0=\epsilon(n,\rho,v)>0$ such that
	if $X$ is the limit space of a sequence of compact
	$n$-manifolds $M_i$ and $\delta_i\to 0$ such that
	$$\op{Ric}_{M_i}\ge -(n-1)\delta_i,\,\op{diam}(M_i)\le 1,\,\op{vol}(M_i)\ge v,\,
	\frac{\op{vol}(B_\rho(x_i^*))}{\svolball{0}{\rho}}\ge 1-\epsilon_0,$$
	then $X$ is isometric to a flat manifold.
\end{lemma}

\begin{proof}
	Arguing by contradiction, assume a sequence $X_i$ such that $X_i$ is not
	isometric to any flat manifold, and $X_i$ is the limit of a sequence of compact
	$n$-manifolds, $M_{ij}\overset{GH}{\longrightarrow}X_i$, as $j\to \infty$, and $M_{ij}$ satisfies the
	conditions in Lemma 3.8 with $\delta_{ij}\to 0$ and $\epsilon_i\to 0$. Passing to a subsequence, we may assume that $X_i\overset{GH}{\longrightarrow}X$, and by a standard diagonal argument
	we may assume a sequence, $M_{ij(i)}\overset{GH}{\longrightarrow}X$. By Theorem 1.2, passing to a subsequence
	we may assume $B_{\frac \rho2}(x_{ij(i)}^*)\overset{GH}{\longrightarrow}\sball{0}{\frac \rho2}$.
	By Corollary 2.2, if $x_{ij(i)}\to x$, then a small ball around $x$ is
	isometric to an Euclidean ball, and thus $X$ is a flat $n$-manifold.
	
	Since $X_i$ is homeomorphic to $M_{ij(i)}$ ((1.7.1)), which, by the
	same reason, is diffeomorphic to $X$, $X_i$ is homeomorphic to $X$.
	Since $\delta_{ij}\to 0$ as $j\to \infty$, $\tilde X_i$ satisfies the
	Splitting property (\cite{CC1}), and thus $\tilde X_i$ is isometric to
	$\Bbb R^{k_i}\times N_i$ and $N_i$ is compact simply connected
	topological manifold. Since $X$ is flat, $\tilde X_i=\Bbb R^n$ and thus
	$X_i$ is flat, a contradiction.
\end{proof}

\begin{proof}[Proof of Theorem 0.4]
	~
	
	Arguing by contradiction, assume $\delta_i\to 0$ and a sequence of compact $n$-manifolds, $M_i\overset{GH}{\longrightarrow}X$, such that $M_i$ is not diffeomorphic to any flat manifold
	and
	$$\op{Ric}_{M_i}\ge -(n-1)\delta_i,\, 1\ge \op{diam}(M_i),\,\volume(M_i)\ge v,\,\frac{\volume
		(B_\rho(x_i^*))}{\svolball0\rho}\ge 1-\epsilon_0,$$
	where $\epsilon(n,\rho,v)$ is from Lemma 3.8. By Lemma 3.8, $X$ is isometric to a flat
	manifold, and by Theorem 1.7 for $i$ large $M_i$ is diffeomorphic to $X$, a contradiction.
\end{proof}

\section{Proof of Theorem D by Assuming Theorem 1.4}

Using Theorem 1.4, we will establish the following result.

\begin{theorem}
	Let $M_i\overset{GH}{\longrightarrow}X$ be a sequence of compact
	$n$-manifolds such that
	$$\op{Ric}_{M_i}\ge -(n-1), \quad
	\op{diam}(M_i)\le d, \quad h(M_i)\geq n-1-\epsilon_i\to n-1.$$
	Then the sequence of Riemannian universal covering
	spaces, $$\begin{CD}(\tilde M_i,\tilde p_i)@>GH>> (\Bbb H^n,o).\end{CD}$$
\end{theorem}

\begin{proof}[Proof of Theorem D by assuming Theorem 4.1]
	~
	
	Arguing by contradiction, assume a sequence of compact $n$-manifolds, $M_i\overset{GH}{\longrightarrow}X$,
	as in Theorem 4.1 such that (3.1.2) holds
	and $M_i$ is not diffeomorphic or not close to any hyperbolic manifold.
	By Theorem 4.1, $\tilde X$ is isomorphic to $\Bbb H^n$. By applying
	Theorem 1.6 on $\tilde M_i$, it is clear that $M_i$
	satisfies the conditions of Theorem B, a contradiction.
\end{proof}

Our proof of Theorem 4.1 is divided into two steps: we first show that $\tilde X$
is isometric to $\Bbb H^k$, $1\le k\le n$ (Lemma 4.4). Then we show that
$\lim_{i\to \infty}h(M_i)=k-1$ (Theorem 4.6), and thus conclude that $k=n$.

To apply Theorem 1.4, we will need to extend an observation in \cite{Li}: if a compact
$n$-manifold of $\op{Ric}_M\ge -(n-1)$ has the maximal volume entropy $n-1$, then
there is a sequence, $r_i\to \infty$, such that for any $\epsilon>0$, (1.5.1) is
satisfied for $L=r_i$ when $i$ large.

\begin{lemma}
	Let $\tilde M$ be a complete Riemannian $n$-manifold such that
	$$h(\tilde M)=\lim_{r\to \infty} \frac{\ln
		\volume(B_r(\tilde{p}))}r\geq n-1-\epsilon.$$
	Then fixing $R>0$ and $\tilde p\in\tilde M$, there exists a
	sequence $r_i\rightarrow \infty$, such that
	$$\lim_{i\rightarrow \infty}\frac{\volume(\partial
		B_{r_i+50R}(\tilde p))}{\volume(\partial B_{r_i-50R}(\tilde
		p))}\geq e^{100R(n-1-\epsilon)},\eqno (4.2.1)$$
	where $e^{100R(n-1)}$ is the limit ratio of the same type in $\Bbb H^n$.
\end{lemma}

\begin{proof}
	Arguing by contradiction, we may assume sufficiently small $\epsilon_0 > 0$ and $r_0 >
	100R$ such that for any $r\geq r_0$,
	$$\frac{\volume(\partial B_{r+50R}(\tilde
		p))}{\volume(\partial B_{r-50R}(\tilde p))}< (1-\epsilon_0)\cdot
	e^{100R(n-1-\epsilon)}.$$
	Then by iteration
	$$\begin{aligned}
	\volume(\partial B_r(\tilde p))&\leq
	(1-\epsilon_0)e^{100R(n-1-\epsilon)}\volume(\partial B_{r-100R}(\tilde
	p))\\
	& \le C(n,r_0,R)\cdot \left((1-\epsilon_0)e^{100R(n-1-\epsilon)}\right)^{\frac{r-r_0}{100R}}.
	\end{aligned}$$
	Thus,
	$$\begin{aligned}
	& h(\tilde{M}) \\
	= &\; \lim_{r\to \infty}\frac{\ln\left(\volume(B_r(\tilde p))\right)}r
	=  \lim_{r\rightarrow \infty}\frac{\ln\left(\int_0^r
		\volume(\partial B_u(\tilde p))du\right)}r\\
	\leq &\; \lim_{r\rightarrow \infty}\frac 1r \ln\left(\int_{r_0}^r
	C(n,r_0,R)\cdot
	\left((1-\epsilon_0)e^{100R(n-1-\epsilon)}\right)^{\frac{u-r_0}{100R}}
	du+\volume(B_{r_0}(\tilde p))\right)\\
	= &\; n-1-\epsilon +\frac{\ln(1-\epsilon_0)}{100R}\\
	< &\; n-1-\epsilon,
	\end{aligned}$$
	a contradiction.
\end{proof}


By Lemma 4.2,  we are able to apply Theorem 1.4 to prove the following:

\begin{lemma}
	Let $M$ be a compact Riemannian $n$-manifold such that $$\op{Ric}_M\ge
	-(n-1), \quad h(M)\geq n-1-\epsilon.$$
	For $R>2\op{diam}(M)=d$, and any $\tilde p\in \tilde M$, there is a connected length metric space $Y$
	such that
	$$d_{GH}\left(B_R\left(\tilde{p}\right),
	B_R\left((0,y)\right)\right)\leq \Psi\left(\epsilon | n,d,
	R\right),$$
	where $B_R((0,y))$ is a metric ball in a warped
	product space $\Bbb R^1\times_{e^s} Y$.
\end{lemma}

\begin{proof}
	By Lemma 4.2, there is $r_i\to \infty$ such that
	$(4.2.1)$ holds.
	Because
	$$\lim_{r\to \infty}\frac{\svolsp{-1}{r+50R}}{\svolsp{-1}{r-50R}}=
	e^{100R(n-1)},$$
	condition (1.5.1) is equivalent to (4.2.1) for $L=r_i> 2R$.
	By Theorem 1.4, for large $i$, $A_{r_i-50R,
		r_i+50R}(\tilde p)$ contains a ball, $B_{2R}(\tilde q)$, such that
	$$d_{GH}\left(B_{2R}\left(\tilde{q}\right),
	B_{2R}\left((0,y)\right)\right)\leq \Psi\left(\epsilon | n,
	R\right),$$
	where $B_{2R}\left((0,y)\right)$ is a metric ball in a warped
	product space $\Bbb R^1\times_{e^s} Y$, and $Y$ is a length metric
	space. Because $R>2\op{diam}(M)$, we may assume that
	$B_{\op{diam}(M)}(\tilde q)$ contains a point $\tilde p'=\gamma(\tilde
	p)$, where $\gamma$ is a deck transformation
	of $\tilde M$. Then $B_R(\tilde p')\subset B_{2R}(\tilde q)$, and
	this completes the proof.
\end{proof}

\begin{lemma}
	Let the assumptions be as in Theorem 4.1.
	Then by passing to a subsequence, $$\begin{CD}(\tilde M_i,\tilde p_i)@>GH>> (\Bbb H^k,o)$ $(k\le n). \end{CD}$$
\end{lemma}

\begin{remark}
	Observe that in Lemma 4.4, if $M_i=M$, then $\tilde M=\Bbb H^n$, and thus $M$ is a hyperbolic manifold.
	This gives a different proof of Theorem 0.3, which does not rely on \cite{LiW} (cf. \cite{LW1},
	\cite{Li}).
\end{remark}

\begin{proof}[Proof of Lemma 4.4]
	Passing to a subsequence, assume that (3.1.2) holds. Fixing any $R>2d$, by Lemma 4.3,
	$$\begin{aligned}
	d_{GH}\left(B_R\left(\tilde p\right), B_{R}\left((0,y_i)\right)\right)
	& \leq
	d_{GH}\left(B_R\left(\tilde p\right),B_R\left(\tilde p_i\right)\right)\\& \hskip4mm +
	d_{GH}\left(B_R\left(\tilde p_i\right),
	B_{R}\left((0,y_i)\right)\right)\\
	&\leq
	\Psi\left(\epsilon_i | n, d, R\right),\end{aligned}$$
	where $B_R((0,y_i))$ is a metric ball in a warped
	product space $\Bbb R^1\times_{e^s} Y_i$.
	Note that $$\begin{CD}(\Bbb R^1\times_{e^s} Y_i, (0,y_i))@>GH>>(\Bbb
	R^1\times_{e^s}Y, (0,y)).\end{CD}$$ Since $R$ is arbitrary, we conclude that $(\tilde X,\tilde p)$
	is isometric to $(\Bbb R^1\times_{e^s}Y,(0,y))$.
	
	Since $\tilde X$ is
	a limit of manifolds of Ricci curvature bounded below, regular points in $\tilde X$ are dense;
	a point is regular if the tangent is unique and isometric to $\Bbb R^k$ for some $k\le n$.
	Without loss of generality, we may assume that $\tilde p$ is a regular point, and thus
	$\lim_{t\to\infty} (e^tY, y)=(\Bbb R^{k-1}, 0)$.
	Via reparametrization of $s'=s-t$,
	$$\begin{aligned} \lim_{t\to \infty}
	\left(\Bbb R^1\times_{e^s}Y,
	(t,y)\right)&=\lim_{t\to\infty}\left(\Bbb R^1\times_{e^{s'}}e^tY,
	(0,y)\right)\\&=\left(\Bbb R^1\times_{e^s}\Bbb R^{k-1}, o\right)=(\Bbb
	H^k,o).\end{aligned}$$
	Since $\tilde X/G$ is compact, for any $t\in
	\Bbb R^1$, there is $\gamma_t\in G$
	such that $d(\gamma_t(\tilde p), (t,y))\le
	\op{diam}(X)\le d$.
	$$(\tilde X, \tilde p)=\lim_{t\to\infty}(\tilde X, \gamma_t(\tilde p))=(\Bbb H^k,o).$$
\end{proof}

\begin{theorem}
	Let $M_i\overset{GH}{\longrightarrow}X$ be a sequence satisfying
	$$\op{Ric}_{M_i}\ge -(n-1), \quad \op{diam}(M_i)\le d,$$
	and the following commutative diagram,
	$$\begin{CD} (\tilde M_i,\tilde p_i,\Gamma_i)@>GH>>(\Bbb H^k,\tilde p,G)
	\\ @VV \pi_i V   @VV \pi  V \\
	(M_i,p_i)@>GH>> (X,p)\end{CD}$$
	Then $\lim_{i\to \infty}h(M_i)=k-1$.
\end{theorem}

Note that Theorem 4.1 follows from Lemma 4.4 and Theorem 4.6.

By Section 1.b, the commutative diagram in Theorem 4.6 yields
the following commutative diagram:
$$\begin{CD} (\tilde M_i,\tilde p_i,\Gamma_i)
@> GH>> (\Bbb H^k,\tilde p,G)\\
@ VV \hat \pi_i V@VV \hat\pi  V\\
(\hat M_i,\hat p_i,\hat \Gamma_i)
@> GH >> (\hat X,\hat p,\hat G)\\
@ VV\bar \pi_i V@VV \bar\pi  V\\
(M_i,p_i) @> GH >> (X,p), \end{CD}$$
where $\Gamma_i\cong\pi_1(M_i)$, $\hat M_i=\tilde
M_i/\Gamma_i^\epsilon$, $\hat X=\Bbb H^k/G_0$,
and $\hat \Gamma_i=\Gamma_i/\Gamma_i^\epsilon\cong G/G_0=\hat G$.
By Lemma 1.13, we may assume an isomorphism $\hat \phi_i: \hat \Gamma_i\to
\hat G$ such that $(\hat h_i,\hat \phi_i,\hat \phi_i^{-1})$ is an $\epsilon_i$-equivariant
GHA on $(B_R(\hat p_i),\hat \Gamma_i(R))$,  $\frac 1{\epsilon_i}>R$. As seen in the proof of Theorem B, $G_0$ is
nilpotent (Theorem 1.9) and thus $G_0=e$ (Theorem 2.5), and thus
$\hat G=G/G_0=G$ is discrete.

\begin{lemma}
	Let the assumptions be as in Theorem 4.6.
	Then for $i$ large, there is a map
	$\hat f_i: (\hat M_i,\hat p_i)\to (\hat X,\hat p)$ such that
	
	\noindent (4.7.1) $\hat f_i$ is an $\epsilon_i$-conjugate, i.e., $|\hat
	f_i(\hat \gamma_i(\hat x_i))\hat \phi_i(\hat \gamma_i)(\hat f_i(\hat x_i))|\leq \epsilon_i$,
	$\hat x_i\in \hat M_i$, $\hat \gamma_i\in \hat \Gamma_i$;
	
	\noindent (4.7.2) for any $R>0$, $\left.\hat f_i\right|_{B_R(\hat
		p_i)}:B_R(\hat p_i)\to B_{\left(1+\frac{\epsilon_i}{60d}\right)R}(\hat
	f_i(\hat p_i))$ is an $\frac{R}{10d}\epsilon_i$-GHA.
\end{lemma}

\begin{proof}
	We first construct a map $\hat f_i: \hat M_i\to \hat X$.
	
	Fix any $R_0>480d$. Let $\hat h_i:
	(B_{\frac{1}{\epsilon_i}}(\hat
	p_i),\hat p_i)\to (\hat X,\hat p)$ be an $\epsilon_i$-equivariant GHA with respect to
	$\hat \phi_i: \hat \Gamma_i\to \hat G$. For $i$ large, we may assume
	that for any $\hat x_i\in B_{R_0}(\hat p_i), \hat \gamma_i\in \hat \Gamma_i(R_0)$,
	$$|\hat h_i(\hat x_i)\hat \phi_i(\hat \gamma_i)^{-1}(\hat h_i(\hat \gamma_i (\hat x_i)))|\leq \epsilon_i.$$
	We now define a map $\hat f_i:\hat
	M_i\to \hat X$ as follows. First, because $\hat \phi_i: \hat \Gamma_i\to \hat G$ is an isomorphism,
	we define $\hat f_i$ on $\hat \Gamma_i(\hat p_i)$ by $\hat f_i(\hat \gamma_i(\hat p_i))=
	\phi_i(\hat \gamma_i)(\hat p)$.
	For any $\hat y_i\in
	\hat M_i\setminus \hat \Gamma_i(\hat  p_i)$, we may assume $\hat\alpha_i\in \hat \Gamma_i$ such
	that $\left|\hat\alpha_i(\hat y_i)\hat p_i\right|\le d$ (note that if $\hat y_i$ is on the boundary of
	a fundamental domain, then $\hat \alpha_i$ is not unique).  We define
	$$\hat f_i(\hat y_i)=\hat \phi_i(\hat\alpha_i)^{-1} (\hat h_i(\hat\alpha_i(\hat y_i))).$$
	If $\hat\beta_i\in \hat\Gamma_i$ satisfies that $\left|\hat\beta_i(\hat y_i)\hat p_i\right|\le d$, then $\hat \beta_i\hat \alpha_i^{-1}\in \hat \Gamma_i(R_0)$, and thus
	$$|\hat h_i(\hat \alpha_i \hat y_i)\hat \phi_i(\hat \alpha_i\hat \beta_i^{-1})(\hat h_i(\hat \beta_i (\hat y_i)))|\leq \epsilon_i.$$
	Since $\hat \phi_i(\alpha_i)^{-1}$ is an isometry, the above implies
	$$|\hat f_i(\hat y_i)\hat \phi_i(\hat \beta_i)^{-1}(\hat h_i(\hat \beta_i (\hat y_i)))|\leq \epsilon_i.$$
	
	(4.7.1) For any $\hat x_i\in \hat M_i, \hat \gamma_i\in \hat \Gamma_i$, let $\hat \alpha_i'$ be the element used
	to define $\hat \gamma_i(\hat x_i)$. Hence,
	$$\hat f_i(\hat \gamma_i(\hat x_i))=\hat \phi_i(\hat \alpha_i')^{-1}\hat h_i(\hat \alpha_i'\hat \gamma_i (\hat x_i))=\hat \phi_i(\hat \gamma_i)\hat \phi_i(\hat \alpha_i'\hat \gamma_i)^{-1}\hat h_i(\hat \alpha_i' \hat \gamma_i (\hat x_i)).$$
	If $\hat \alpha_i$ denotes the element defining $\hat x_i$, then we may view $\hat \alpha_i'\hat \gamma_i$
	as $\hat \beta_i$ as in the above discussion. By now we can conclude that
	$$|\hat f_i(\hat \gamma_i(\hat x_i))\hat \phi_i(\hat \gamma_i)(\hat f_i(\hat x_i))|\leq \epsilon_i.$$
	
	$(4.7.2)$. Since $\hat f_i$ is $\epsilon_i$-onto
	from $B_{R_0}(\hat p_i)$ to $B_{R_0}(\hat f_i(\hat p_i))$ and
	$\hat f_i$ is $\epsilon_i$-conjugate, $\hat f_i$ is $2\epsilon_i$-onto (For any $\hat x\in \hat X$, there is $\hat\gamma\in \hat G$, such that $\hat \gamma(\hat x)\in B_{R_0}(\hat p)$. Then there is $\hat \gamma_i\in \hat \Gamma_i, \hat x_i\in B_{R_0}(\hat p_i)$, such that $\hat \phi_i(\hat\gamma_i)=\hat \gamma$ and $|\hat f_i(\hat x_i)\hat \gamma(\hat x)|\leq \epsilon_i$. Since $\hat f_i$ is $\epsilon_i$-conjugate,
	$|\hat f_i(\hat \gamma_i^{-1}(\hat x_i))\hat x|\leq |\hat \phi_i(\hat\gamma_i^{-1})\hat f_i(\hat x_i)\hat x|$ $+\epsilon_i=|\hat \gamma^{-1}\hat f_i(\hat x_i)\hat x|+\epsilon_i\leq 2\epsilon_i$).
	
	For any $R>R_0$ and any $\hat x_i,
	\hat y_i \in B_R(\hat p_i)$,
	we shall estimate
	$$\left||\hat x_i\hat y_i|-|\hat f_i(\hat x_i)\hat f_i(\hat y_i)|\right|.$$
	Let $\hat c:[0,l]\to \hat
	M_i$ $(l=|\hat x_i\hat y_i|)$ be a minimal geodesic
	connecting $\hat x_i$ and $\hat y_i$
	parametrized by arc length, and let
	$0=t_0<t_1<\cdots<t_s=l$ of $[0,l]$ be a partition such that
	$t_{j+1}-t_j=\frac{R_0}{2}$ $(0\le j<s-1)$ and
	$t_{s}-t_{s-1}\le\frac{R_0}{2}$. Then
	$s\le \frac{2l}{R_0}$
	and $|\hat c(t_j)\hat c(t_{j+1})|\le \frac{R_0}{2}$. For each $j$, there is
	$\hat \gamma_j\in \hat \Gamma_i$ such that
	$B_{R_0}(\hat\gamma_j(\hat p_i))$ contains
	$\hat c|_{[t_j,t_{j+1}]}$.
	Because $\hat f_i$ is an $\epsilon_i$-conjugate and $\epsilon_i$-GHA on
	$B_{R_0}(\hat p_i)$, we derive
	$$\begin{aligned}
	& \left|\left|\hat c(t_j) \hat c(t_{j+1})\right|-\left|\hat f_i(\hat c(t_j))\hat f_i(\hat c(t_{j+1}))\right|\right|\\
	= &\left|\left|\hat \gamma_j^{-1}(\hat c(t_j))
	\hat\gamma_j^{-1}(\hat c(t_{j+1}))\right| -
	\left|\hat\phi_i(\hat\gamma_j^{-1})\hat f_i
	(\hat c(t_j)) \hat\phi_i(\hat\gamma_j^{-1}) \hat
	f_i(c(t_{j+1}))\right|\right|\\
	\leq &\left|\left|\hat \gamma_j^{-1}(\hat c(t_j))
	\hat\gamma_j^{-1}(\hat c(t_{j+1}))\right|-\left|\hat f_i
	(\hat\gamma_j^{-1}\hat c(t_j)) \hat
	f_i(\hat\gamma_j^{-1}\hat c(t_{j+1}))\right|\right|+2\epsilon_i\\
	\le & 3\epsilon_i.
	\end{aligned}$$
	Then
	$$\left|\hat f_i(\hat x_i)\hat f_i(\hat y_i)\right|\le \sum_{j} \left|\hat f_i
	(\hat c(t_j))\hat f_i(\hat c(t_{j+1}))\right|\le
	\left(1+\frac{6}{R_0}\epsilon_i\right)\left|\hat x_i\hat
	y_i\right|.$$
	
	To establish the opposite inequality, note that a minimal geodesic
	between $\hat f_i(\hat x_i)$ and $\hat f_i(\hat y_i)$ may not lie
	in the image of $\hat f_i$. Since $\hat f_i$ is $2\epsilon_i$-onto,
	we may replace the partition points by nearby points in $\hat f_i(\hat M_i)$.
	Similar to the above estimate we derive
	$$\left|\hat x_i\hat y_i\right|\le
	\left(1+\frac{24}{R_0}\epsilon_i\right)\left|\hat f_i(\hat x_i)\hat
	f_i(\hat y_i)\right|.$$
	Now $(4.7.2)$ follows by taking $R_0=480d$.
\end{proof}

Let $\pi: (\tilde M,\tilde p)\to (M,p)$ be the Riemanian covering space, and let $\Gamma=\pi_1(M,p)$.
Observe that if $\op{diam}(M)\le d$, then for any $R>0$,
$$\frac{\volume (B_{R-d}(\tilde
	p))}{\volume
	(B_{d}(p))}\le |\Gamma(R)|\le \frac{\volume
	(B_{R+d}(\tilde p))} {\volume (B_{d}(p))},$$
and thus
$$h(M)=\lim_{R\to \infty}\frac{\ln \op{vol}(B_R(\tilde p))}R=\lim_{R\to \infty}\frac{\ln |\Gamma(R)|}{R}.$$

\begin{proof}[Proof of Theorem 4.6]
	~
	
	Let $\epsilon>0$ satisfy that $\Gamma_i^\epsilon \overset{GH}{\longrightarrow} G_0$ (see Lemma 1.13). By Theorem
	2.5, $G_0=e$.
	Then $\Gamma_i^\epsilon(\tilde p_i)\to \tilde p$, and thus $\Gamma_i^\epsilon$ is finite
	when $i$ large. For $\gamma_i\in \Gamma_i(R)$, we may assume
	$\gamma_i\in \alpha_i\Gamma_i^\epsilon$.
	Observe that $\alpha_i$ can be chosen
	so that $\hat \alpha_i\in \hat \Gamma_i(R)$,
	where $\hat \alpha_i$ denotes the projection of $\alpha_i$ in $\hat \Gamma_i$.
	Assume that $|\Gamma_i^\epsilon|=C_i$. Then
	$$|\hat\Gamma_i(R)|\le |\Gamma_i(R)|\le |\hat\Gamma_i(R)|\cdot|\Gamma_i^\epsilon|\le C_i\cdot |\hat \Gamma_i(R)|. \eqno (4.6.1)$$
	
	We claim that
	$$C_1 e^{\left(k-1\right)\left(1-
		\frac{\epsilon_i}{10d}\right)R} \le |\hat \Gamma_i(R)|\le C_2 e^{\left(k-1\right)\left(1+
		\frac{\epsilon_i}{10d}\right)R}\eqno (4.6.2)$$
	Combining (4.6.1) and (4.6.2), we derive
	$$\left|\frac{1}{k-1}\cdot
	h(M_i)-1\right|=\left|\frac{1}{k-1}\cdot\lim_{R\to \infty}\frac {\ln
		|\Gamma_i(R)|}R -1\right|\le \frac{\epsilon_i}{10d}.$$
	
	We now verify (4.6.2).
	Let $\hat f_i$ be in Lemma 4.7, $\hat p=\hat f_i(\hat p_i)$. By (4.7.2) for any $R>0$,
	$$\hat G(\hat p)\cap B_{\left(1-\frac{\epsilon_i}{10d}
		\right)R}(\hat p)\subset \hat
	f_i\left(\hat\Gamma_i(R)(\hat p_i)\right)\subset
	\hat G(\hat p)\cap B_{\left(1+\frac{\epsilon_i}{10d}
		\right)R}(\hat p).$$
	
	Without loss of generality, we may assume $\delta>0$ such that $\hat G$ has a trivial isotropy
	group in $B_{2\delta}(\hat p)$ and $\hat G(\hat p)\cap B_{2\delta}(\hat p)=\{\hat p\}$.
	Together with the fact that $\hat f_i$ is $\epsilon_i$-conjugate, we have that
	$$|\hat G(\hat p)\cap B_{\left(1-\frac{\epsilon_i}{10d}\right)R}(\hat p)|
	\le|\hat \Gamma_i(R)|
	\le |\hat G(\hat p)\cap B_{\left(1+\frac{\epsilon_i}{10d}\right)R}(\hat p)|
	\eqno (4.6.3)$$
	Counting points in $\hat G(\hat p)\cap B_R(\hat p)$, we get
	$$
	\frac{\svolball{-1}{R}} {\svolball{-1}{d}}\le
	|\hat G(\hat p)\cap B_R(\hat p)|\le
	\frac{\svolball{-1}{R}} {\svolball{-1}{\delta}},\quad
	B_R(\hat p)=\sball{-1}{R}.\eqno (4.6.4)$$
	By now, (4.6.2) follows from (4.6.3) and (4.6.4).
\end{proof}

\begin{proof}[Proof of Theorem 0.5]
	~
	
	The proof is similar to the proof of Theorem 4.6,
	because $\dim(M)=n$. Hence, we will only briefly describe the proof.
	
	First, since $\dim(M)=n$, $G_0=e$, and since $\Gamma_i^\epsilon\overset{GH}{\longrightarrow}e$,
	by Lemma 2.3 we conclude that for $i$ large, $\Gamma_i^\epsilon=e$.
	By Lemma 1.13, we see that $\hat \Gamma_i=\Gamma_i/\Gamma_i^\epsilon
	\cong G/G_0=G$. Assume $(h_i, \phi_i, \phi^{-1})$ be $\epsilon_i$-equivariant
	GHA with $\epsilon_i\to 0$, where $\phi_i:\Gamma_i\to G$ is an isomorphism.
	
	Following the proof of Lemma 4.7 with $\hat M_i=\tilde M_i$ and $\hat X=\tilde X=\tilde M$,
	via the center of mass method we construct a map, $\tilde f_i:
	(\tilde M_i,\tilde p_i,\Gamma_i)\to (\tilde M,\tilde p,G)$, such that
	(4.7.1) and (4.7.2) hold. By the estimate for $\hat \Gamma_i$ in
	the proof of Theorem 4.6, we get the desired result.
\end{proof}

\begin{proof}[Corollary 0.6]
	~
	
	(0.6.1) $\Rightarrow$ (0.6.3): By Theorem 0.5.
	
	(0.6.3) $\Rightarrow$ (0.6.2): By Theorem 4.1, $\tilde M$ is close to $\Bbb H^n$. By
	Theorem 1.6 we see that (0.6.2) is satisfied.
	
	(0.6.2) $\Rightarrow$ (0.6.1): By Theorem B.
\end{proof}

\section{Proof of Theorem 1.4}

Our approach to Theorem 1.4 is based on the following functional criterion for warped product metric
by Cheeger-Colding (see Theorem 5.1).

Let $N$ be a Riemannian $(n-1)$-manifold, let $k:(a,b)\to \Bbb R$ be a smooth positive function,
and let $(a,b)\times_k N$ be the $k$-warped product whose Riemannian tensor is
$$g = dr^2 + k^2(r)g_N.$$
Then the function, $f = -\int_r^b k(u)du,$
satisfies $$\op{Hess}f = k'(r)g.$$
Conversely, let $(M, g)$ be a Riemannian manifold and let $r:M\to \Bbb R$ be the distance function to
a compact subset of $M$.  If there is a smooth function $f: M\to \Bbb R$ satisfying
$\nabla f \ne 0$ and
$$\op{Hess}f = h\cdot g$$ on $A_{a,b}=r^{-1}((a,b))$, where  $h:M\to \Bbb R$ is a smooth function,
then $f$ is constant on each level set of $r$ between $a$
and $b$, and the Riemannian metric tensor in the annulus
$A_{c,d}$ ($a< c< d< b$) is a warped product (cf. \cite{CC1}),
$$g = dr^2 + (f'(r))^2 \tilde{g}.$$

Cheeger-Colding proved that if $\op{Hess}f = k'(r) g$
holds approximately ``in the $L_1$-sense'', then the warped product structure
of $A_{c,d}$ almost holds ``in the Gromov-Hausdorff sense'' \cite{CC1}.

\begin{theorem}[\cite{CC1}]
	Let $M$ be a Riemannian manifold with $\op{Ric}_M\geq -(n-1)H$, let $r$ be a distance
	function to a compact subset in $M$, let $k:\Bbb R\to \Bbb R$ be a positive smooth function
	and let $f=-\int_r^bk(u)du$. For $0<\alpha'<\alpha$, let $A_{a+\alpha,b-\alpha}\subset
	A_{a+\alpha',b-\alpha'}$ be two annuluses with respect to $r$. Let $d^{\alpha'}$ be the intrinsic
	metric in $A_{a+\alpha', b-\alpha'}$, and let $d^{\alpha',\alpha}=d^{\alpha'}|_{A_{a+\alpha, b-\alpha}}$.
	Assume
	
	\noindent (5.1.1) for the metric $d^{\alpha',\alpha}$, $\op{diam}(A_{a+\alpha, b-\alpha})\leq D$,
	
	\noindent (5.1.2) for $0<\delta<\alpha'$ and all $x\in r^{-1}(a+\alpha')$, there exists $y\in r^{-1}(b-\alpha')$ such that
	the intrinsic distance between $x$ and $y$ in $A_{a+\alpha'-\delta,
		b-\alpha'+\delta}$ satisfies $$d^{\alpha'-\delta}(x,y)\leq b- a - 2\alpha' + \delta.$$
	
	\noindent (5.1.3) there is $\tilde{f} : A_{a,b}\rightarrow \Bbb R$ satisfying
	
	(5.1.3.1) $|\tilde{f} - f |<\delta$ for all $x\in A_{a+\alpha', b-\alpha'}$,
	
	(5.1.3.2) $-\kern-1em\int_{A_{a, b}}|\nabla \tilde{f}-\nabla f|\le \delta$,
	
	(5.1.3.3) $-\kern-1em\int_{A_{a+\alpha', b-\alpha'}}|\op{Hess}\tilde{f}-
	k'(r)g|\le \delta$,
	
	\noindent Then there exits a metric space $X$, with $\op{diam}(X)\le C(a, b,
	\alpha, \alpha', f, D, H)$, such that for the restricted metric $d^{\alpha, \alpha'}$ on $A_{a+\alpha,b-
		\alpha}$,
	$$d_{GH}(A_{a + \alpha, b - \alpha}, (a + \alpha, b - \alpha)\times_{k}X)\le
	\Psi(\delta | a, b, \alpha, \alpha', n, f, D, H).$$
\end{theorem}

We will only present a proof of Theorem 1.4 for $H<0$, because
a proof for $H=0$ follows the same argument with a minor modification.
By a rescaling, without loss of generality we assume $H=-1$.

From the proof of Theorem 5.1 (see Proposition 2.80 and Theorem 3.6 in \cite{CC1}),
we observe the following: If (5.1.2) holds on $B_\rho(q)\subset
A_{a+\alpha, b-\alpha}(p)$, and one can find $\tilde f$ such that (5.1.3)
holds, then $$d_{GH}(B_\rho(p),B_\rho(0,y))<\Psi(\delta |\rho, n, f, H),$$
where $B_\rho(0,y)\subset (a + \alpha, b - \alpha)\times_{k}X$
for some metric space $X$.

In view of Theorem 1.4, we choose $f=e^u, u(x)=|xp|-|pq|$, for some
$q\in A_{L-R,L+R}(p)$ such that $B_\rho(q)$ satisfies (5.1.2),
and $\tilde f$ is the solution of
$$\begin{cases}\Delta \tilde f = ne^u, & \text{in } B_\rho(q);\\
\tilde f =f, &\text{on } \partial B_\rho(q).\end{cases}\eqno (5.2)$$

Our strategy is to select balls in $A_{L-2R,L+2R}(p)$ such that
(5.1.2) holds on each ball (see Lemmas 5.4 and 5.5), which also
satisfies an additional property (see Lemma 5.8) so that we are able to verify (5.1.3)
(see Lemma 5.9).

From the above discussion, the following theorem implies Theorem 1.4.

\addtocounter{theorem}{1} 
\begin{theorem}
	Let the assumptions be as in Theorem 1.4. Given $0<\alpha<1$, there are disjoint metric balls,
	$B_{\rho}(q_i)\subset A_{L-R, L+R}(p)$, satisfying (1.4.2) and the following:
	
	\noindent (5.3.1) for $x\in B_\rho(q_i)$, there is $y\in \partial B_{L+ R}(p)$
	satisfying $|p x|+ |x y|\leq |p y|+ \Psi(\epsilon, L^{-1} |n,\rho,R)$;
	
	\noindent (5.3.2) for each $q_i$, let $u(x)=|x p|-|q_i p|$, there is a
	smooth function $\tilde f$ satisfying
	
	(5.3.2.1) $|\tilde f-e^u|<\Psi(\epsilon, L^{-1} |n, R,\rho)$ for all $x\in B_{(1-\alpha)\rho}(q_i)$.
	
	(5.3.2.2) $-\kern-1em\int_{B_\rho(q_i)}|\nabla \tilde f-\nabla e^u|^2\le
	\Psi(\epsilon, L^{-1} |n,R,\rho)$.
	
	(5.3.2.3) $-\kern-1em\int_{B_{(1-\alpha)\rho}(q_i)}|\op{Hess}\tilde f-e^u|^2
	\le \Psi(\epsilon, L^{-1}|n,R,\alpha,\rho)$.
\end{theorem}

From now on, without mention explicitly we always assume the condition (1.5.1) and
denote $\epsilon=\Psi(\epsilon|n,H,R)$.

Let $E$ be a maximal subset of $\{q_i, \,B_{\rho}(q_i)\subset A_{L-R, L+R}(p)\}$
such that for all $q_{i_1}\neq q_{i_2} \in E$, $B_{\rho}(q_{i_1})\cap B_{\rho}(q_{i_2})=\emptyset$.
Let $F=\bigcup_{q_i\in E}B_{\rho}(q_i)$. We shall choose $q_i\in E$ such that (5.3.1) and (5.3.2)
hold on $B_\rho(q_i)$.

\begin{lemma}
	For $L$ sufficiently large,
	$$\frac{\volume(F)}{\volume(A_{L-R,L+R}(p))}\ge (1-\Psi(\epsilon,L^{-1}|n,\rho,R))e^{-(n-1)\rho}\cdot \frac{\svolball{-1}{\rho}}{\svolball{-1}{2\rho}}.$$
\end{lemma}

\begin{proof}
	Let $G=\bigcup_{q_i\in E}B_{2\rho}(q_i)$. By the maximality of $E$, we have that,
	$$A_{L-R+\rho, L+R-\rho}(p) \subset G.$$ For $L-R<r<L+R$, by (1.5.1) and
	the Bishop-Gromov relative volume comparison, we get
	$$\frac{\volume(\partial B_r(p))}{\svolsp{-1}{r}}\geq (1-\epsilon)\frac{\volume
		(\partial B_{L-R}(p))}{\svolsp{-1}{L-R}}.$$
	Plugging the above into the integrant in the following quotient, together
	with the Bishop-Gromov relative volume comparison, we derive
	$$\begin{aligned} \frac{\volume(G)}{\volume(A_{L-R, L+R}(p))}\ge &
	\frac{\volume(A_{L-R+\rho, L+R-\rho}(p))}{\volume(A_{L-R, L+R}(p))}\\ = & \frac{\int_{L-R+\rho}^{L+R-
			\rho}\volume(\partial B_r(p))dr}{\int_{L-R}^{L+R}\volume(\partial B_r(p))dr}\\
	\geq & \frac{\int_{L-R+\rho}^{L+R-\rho}(1-\epsilon)\frac{\volume(\partial B_{L-R}(p))}{\svolsp{-1}
			{L-R}}\svolsp{-1}{r}dr}{\int_{L-R}^{L+R}\frac{\volume(\partial B_{L-R}(p))}{\svolsp{-1}{L-R}}\svolsp{-1}{r}dr} \\
	= & (1-\epsilon)\frac{\svolann{-1}{L-R+\rho, L+R-\rho}}{\svolann{-1}{L-R, L+R}}\\
	\geq & (1-\Psi(\epsilon, L^{-1}|n,\rho,R))e^{-(n-1)\rho}.
	\end{aligned} \eqno (5.4.1)$$
	Again applying Bishop-Gromov relative volume comparison to the numerator of
	the quotient,
	$$\frac{\volume(F)}{\volume(G)}\geq \frac{\sum_{q_i\in E} \volume(B_{\rho}(q_i))}{\sum_{q_i\in E}
		\volume(B_{2\rho}(q_i))}\geq \frac{\svolball{-1}{\rho}}{\svolball{-1}{2\rho}}.
	\eqno (5.4.2)$$
	The desired result follows from (5.4.1) and (5.4.2).
\end{proof}

Next, we show that the balls in $F$ on which (5.3.1) and (5.3.2)
hold have a total volume almost equals to the volume of $F$.

Let $S\subset A_{L-R,L+R}(p)$ consist of interior points of minimal
geodesics $c_y$ from $p$ to $y\in \partial B_{L+R}(p)$, i.e.,
$$S=\{x\in A_{L-R,L+R}(p)\cap c_y,\, y\in \partial B_{L+R}(p)\}.$$
Fixing $0<\eta<1$ (which will be specified later), let
$$E'(\eta)=\left\{q_i\in E,\,\,\, \;\frac{\volume(B_\rho(q_i)\setminus S)}{\volume(B_\rho(q_i))}<\eta\right\},$$
and let $F'(\eta)=\bigcup_{q_i\in E'(\eta)}B_{\rho}(q_i)$.

\begin{lemma}
	Let $F'(\eta)$ be defined in the above. Then
	$$\frac{\volume(F'(\eta))}{\volume(F)}\ge1-\eta^{-1}\Psi_1(\epsilon|n,R,\rho).$$
\end{lemma}

\begin{proof}
	Since for any $q_i\in E\setminus E'(\eta)$,
	$$\frac{\volume(B_\rho(q_i)\setminus S)}{\volume(B_\rho(q_i))}\geq \eta,$$
	adding $\op{vol}(B_\rho(q_i))$ over $q_i$'s in $E\setminus E'(\eta)$ we derive
	$$\frac{\volume((F\setminus F'(\eta))\setminus S)}{\volume(F\setminus F'(\eta))}\geq \eta.\eqno (5.5.1)$$
	By Bishop-Gromov relative volume comparison and (1.5.1),
	$$\begin{aligned}
	\frac{\volume(S)}{\svolann{-1}{L-R, L+R}} &\overset{(BG)}{\ge} \frac{\volume(\partial B_{L+R}(p))}{\svolsp{-1}{L+R}}\\
	&\overset{(1.5.1)}{\ge} (1-\epsilon)\frac{\volume(\partial B_{L-R}(p))}{\svolsp{-1}{L-R}}\\
	&\overset{(BG)}{\ge} (1-\epsilon) \frac{\volume(A_{L-R, L+R}(p))}{\svolann{-1}{L-R, L+R}}
	\end{aligned}$$
	By (5.5.1) and Lemma 5.4,
	$$\begin{aligned} &\frac{\volume(S)}{\volume(A_{L-R, L+R}(p))}\\
	=&\; 1-\frac{\volume(A_{L-R,L+R}(p)\setminus S)}
	{\volume(A_{L-R, L+R}(p))}\\
	\leq &\; 1-\frac{\volume((F\setminus F'(\eta))
		\setminus S)}{\volume(F\setminus F'(\eta))}\frac{\volume(F\setminus F'(\eta))}{\volume(F)}\frac{\volume(F)}{\volume(A_{L-R, L+R}(p))}\\
	\leq &\; 1-\eta \cdot c(n,R,\rho)\cdot \frac{\volume(F\setminus F'(\eta))}{\volume(F)}.
	\end{aligned}$$

	Combining the two estimates on $\volume(S)$, we derive
	$$\frac{\volume(F'(\eta))}{\volume(F)}\geq 1-\eta^{-1}\cdot \epsilon\cdot c^{-1}(n,R,\rho).$$
\end{proof}

For $p\in M$, let $C_p$ be the cut locus of $p$.
\begin{lemma}
	Let the assumptions be as in Theorem 1.4,
	and let $r(x)=|px|$. Then for a.e. $L$,
	$$-\kern-1em\int_{A_{L-R,L+R}(p)\setminus C_p}\left|\Delta r-(n-1)\right|
	\le \Psi_2(\epsilon, L^{-1}|n, R).$$
\end{lemma}

\begin{proof}

As the proof of \cite[Theorem 4.1]{Ch}, for each $\delta>0$, take a smooth neighborhood $U_{\delta}$ of $C_p$, such that for $x\in A_{L-R, L+R}(p)\cap \partial U_{\delta}$, $\nabla r$ points into $U_{\delta}$. 
Then
\begin{eqnarray*}
& & \int_{A_{L-R, L+R}(p)\setminus C_p} \Delta r\\
& = &\lim_{\delta\to 0}\int_{A_{L-R, L+R}(p)\setminus U_{\delta}} \Delta r\\
&=&  \lim_{\delta\to 0}\int_{\partial B_{L+R}(p)\setminus U_{\delta}}\left<\nabla r, \nabla r\right>+\lim_{\delta\to 0}\int_{\partial B_{L-R}(p)\setminus U_{\delta}}\left<\nabla r, -\nabla r\right> \\
& & +\lim_{\delta\to 0}\int_{A_{L-R, L+R}(p)\cap \partial U_{\delta}}\left<\nabla r, N\right>\\
&=& \volume(\partial B_{L+R}(p)\setminus C_p)-\volume(\partial B_{L-R}(p)\setminus C_p)+\lim_{\delta\to 0}\int_{A_{L-R, L+R}(p)\cap \partial U_{\delta}}\left<\nabla r, N\right>\\
&\geq &\volume(\partial B_{L+R}(p)\setminus C_p)-\volume(\partial B_{L-R}(p)\setminus C_p),
\end{eqnarray*}
where $N$ is the inward normal to $\partial U_{\delta}$, and thus $\left<\nabla r, N\right> >0$.
Note that $\volume(C_p)=0$. By $$\volume(B_{L+R}(p))=\int_0^{L+R} \volume(\partial B_r(p))dr =\int_0^{L+R}\volume(\partial B_r(p)\setminus C_p)dr,$$ we have that for a.e $r>0$, $\volume(\partial B_r(p))=\volume(\partial B_r(p)\setminus C_p)$. 
And thus for a.e. $L$,
$$\int_{A_{L-R, L+R}(p)\setminus C_p} \Delta r \geq \volume(\partial B_{L+R}(p))-\volume(\partial B_{L-R}(p)).$$

	Then by (1.5.1)
	and
	$$\lim_{L\rightarrow \infty}\frac{\svolsp{-1}{L-R}}{\svolann{-1}{L-R, L+R}}= \frac{n-1}{e^{2R(n-1)}-1},$$
	we derive
	$$\begin{aligned}-\kern-1em\int_{A_{L-R, L+R}(p)\setminus C_p}\Delta r = & \; \frac{\volume(\partial B_{L+R}(p))- \volume(\partial B_{L-R}(p))}
	{\volume(A_{L-R, L+R}(p))}\\ \geq & \; \left((1-\epsilon)\frac{\svolsp{-1}{L+R}}{\svolsp{-1}{L-R}}-1 \right)\frac{\volume
		(\partial B_{L-R}(p))}{\volume(A_{L-R, L+R}(p))}\\ \geq & \left((1-\epsilon)\frac{\svolsp{-1}{L+R}}{\svolsp{-1}
		{L-R}}-1 \right)\frac{\svolsp{-1}{L-R}}{\svolann{-1}{L-R, L+R}}\\
	\geq &\; (1-\Psi(\epsilon, L^{-1}|n, R))(n-1).
	\end{aligned}
	$$
	Let $\underline{\Delta}$ denote the Laplacian on $\Bbb H^n$. By Laplace comparison, we derive
	$$\begin{aligned}
	-\kern-1em\int_{A_{L-R, L+R}(p)}\Delta r &\leq -\kern-1em\int_{A_{L-R, L+R}(p)}\underline{\Delta} r\\
	& = -\kern-1em\int_{A_{L-R, L+R}(p)}(n-1)\frac{\cosh r}{\sinh r}\\
	& \le (1+\Psi( L^{-1}|n,R))(n-1).
	\end{aligned} \eqno (5.6.1)$$
	The desired estimate then follows from the above two estimates for $-\kern-1em\int_{A_{L-R, L+R}(p)\setminus C_p}\Delta r$.
\end{proof}
In the following Lemma 5.7 and Lemma 5.8, we will use relative volume comparison to derive the estimate of $\Delta r$ in small balls. And we will write $B_R(x)$ instead of $B_R(x)\setminus C_p$ for simplicity.
\begin{lemma}
	$$-\kern-1em\int_{F}|\Delta r-(n-1)|\le \Psi_3(\epsilon, L^{-1}|n,R,\rho).$$
\end{lemma}
\begin{proof}
	By Lemma 5.4 and Lemma 5.6, we have that
	$$\begin{aligned}
	-\kern-1em\int_{F}\Delta r = &\; \frac{\volume(A_{L-R, L+R}(p))}{\volume(F)}\left(-\kern-1em\int_{A_{L-R, L+R}(p)}
	\Delta r\right)-\frac{\int_{A_{L-R, L+R}(p)\setminus F} \Delta r}{\volume(F)}\\
	\geq &\; (1-\Psi(\epsilon, L^{-1}|n,R))(n-1)\frac{\volume(A_{L-R, L+R}(p))}{\volume(F)}\\
	& -(n-1 +\Psi(L^{-1}|n, R))\frac{\volume(A_{L-R, L+R}(p)\setminus F)}{\volume(F)}\\
	\geq &\; (1-\Psi(\epsilon, L^{-1}|n,R,\rho))(n-1).
	\end{aligned}$$
	As in (5.6.1), we derive
	$$-\kern-1em\int_{F}\Delta r\leq (1+\Psi(L^{-1}|n,R))(n-1).$$
\end{proof}

Let $$\begin{aligned}
&\Psi(\epsilon,L^{-1}|n,R,\rho)\\
=&\; \max\left\{\Psi_1
(\epsilon|n, R,\rho),\Psi_2(\epsilon,L^{-1}|n, R),\Psi_3
(\epsilon,L^{-1}|n, R,\rho)\right\}.
\end{aligned}$$

\begin{lemma}
	Let $$E''(\eta)=\left\{q_i\in E,\,\, -\kern-1em\int_{B_\rho(q_i)} \left|\Delta r-(n-1)\right|<\eta^{-1}\Psi(\epsilon,L^{-1}|n,R,\rho)\right\},$$
	and let $F''(\eta)=\bigcup_{q_i\in E''(\eta)}B_{\rho}(q_i)$.
	Then
	$$\frac{\volume(F''(\eta))}{\volume (F)}\ge 1-\eta.$$
\end{lemma}

\begin{proof}
	By Lemma 5.7, we derive
	$$\begin{aligned}
	\Psi(\epsilon, L^{-1}|n,R,\rho)\geq &-\kern-1em\int_{F}|\Delta r-(n-1)|\\
	=& \frac{1}{\volume(F)}\left(\sum_{ E''(\eta)}\volume(B_{\rho}(q_i))-\kern-1em\int_{B_{\rho}(q_i)}|\Delta r-(n-1)|\right.\\
	& \left.+\sum_{ E\setminus E''(\eta)}\volume(B_{\rho}(q_i))-\kern-1em\int_{B_{\rho}(q_i)}|\Delta r-(n-1)|\right)\\
	\geq &\frac{1}{\volume(F)}\left(0 + \eta^{-1}\Psi(\epsilon,L^{-1}|n,R,\rho)\volume(F\setminus F''(\eta))\right)\\
	=& \eta^{-1}\Psi(\epsilon,L^{-1}|n,R,\rho)\frac{\volume(F\setminus F''(\eta))}{\volume(F)},
	\end{aligned}$$
	i.e., $$\frac{\volume(F\setminus F''(\eta))}{\volume(F)}\leq \eta.$$
\end{proof}

We now specify $\eta=\sqrt{\Psi(\epsilon,L^{-1}|n,R,\rho)}$.
Then $F'(\eta)\cap F''(\eta)$ satisfies (1.4.2). By Bishop-Gromov
relative volume comparison, (5.3.1) holds on balls in $F'(\eta)$.

To verify (5.3.2) on $B_\rho(q_i)$ for $q_i\in E'(\eta)\cap E''(\eta)$, we will use the standard comparison functions (see \cite{Ch} for more details). Let
$$\underline{U}(r) = \int_0^r sn_H^{1-n}(s)\left(\int_0^s sn_H^{n-1}(u)du\right)ds,$$
$$\underline{G}(r) = \frac{1}{\omega^{n-1}}\int_r^{\infty}sn_H^{1-n}(s)ds,$$
where $\omega^{n-1}=\op{vol}(S_1^{n-1})$. For fixed $d> 0$,
$$\underline{U}_d(r) = \underline{U}(r) - \underline{U}(d),\quad \underline{G}_d(r) = \underline{G}(r) - \underline{G}(d),$$
$$\underline{L}_d(r) = -\frac{\underline{U}'(d)}{\underline{G}'(d)}\underline{G}_d(r) + \underline{U}_d(r).$$
Then $\underline{L}'_d(r)\leq 0, r \in [0,d]$, $\underline{\Delta} \underline{L}_d(r) = 1$, $\underline{\Delta}\underline{U}_d=1$ and $\underline{U}_d'\ge 0$.

\begin{lemma}
	$(5.3.2)$ holds for each $q_i\in E''(\eta)$.
\end{lemma}

\begin{proof}
	For $q=q_i\in E''(\eta)$, let $u(x)=|px|-|pq|$. By Lemma 5.8,
	$$-\kern-1em\int_{B_\rho(q)\setminus C_p} \left|\Delta u-(n-1)\right|<\Psi(\epsilon,L^{-1}|n,R,\rho).$$
	Let $\tilde f$ be the solution of (5.2). Then
	\begin{eqnarray*}
	-\kern-1em\int_{B_{\rho}(q)\setminus C_p}|\Delta (\tilde f-e^u)|&=&-\kern-1em\int_{B_{\rho}(q)\setminus C_p} |n e^u - e^u(|\nabla u|^2 + \Delta u)|\\
	&=& -\kern-1em\int_{B_{\rho}(q)\setminus C_p}e^u|n-1-\Delta u|\\
	&\leq & \Psi(\epsilon, L^{-1}|n,R,\rho).
	\end{eqnarray*}
	
	By maximum principle, $$\Delta (\tilde f-ne^{-2R}\underline{U}_{4R}(u+2R))\geq 0,$$ and $$\Delta (\tilde f - n e^{2R}\underline{L}_{5R}(u+2R))\leq 0,$$ we have that
	$|\tilde f - e^u| \leq c(n,R,\rho)$. We then derive (5.3.2.2) as follows:
	$$\begin{aligned}
	& -\kern-1em\int_{B_{\rho}(q)}|\nabla \tilde f -\nabla e^u|^2\\
	= & -\kern-1em\int_{B_{\rho}(q)}\left<\nabla \tilde f,\nabla \tilde f -\nabla e^u\right>-\left<\nabla e^u, \nabla \tilde f -\nabla e^u\right>\\
	=  & 	-\kern-1em\int_{B_{\rho}(q)} -\Delta \tilde f (\tilde f-e^u) + -\kern-1em\int_{B_{\rho}(q)\setminus C_p}\Delta e^u(\tilde f- e^u)\\
	&- \lim_{\delta \to 0}\frac{1}{\volume(B_{\rho}(q))}\int_{\partial U_{\delta}\cap B_{\rho}(q)}\langle
	\nabla e^u, v\rangle (\tilde f -  e^u)\\
	\leq &- -\kern-1em\int_{B_{\rho}(q)\setminus C_p}(\Delta f-\Delta e^u)(\tilde f- e^u)\\
	\leq & \Psi(\epsilon, L^{-1}|n,R,\rho),\end{aligned}$$
	where $v$ is the inward normal vector to $\partial U_{\delta}\cap B_{\rho}(q)$, $\left<\nabla u, v\right> >0$, and $ U_{\delta}$ is a $\delta$-tube neighborhood of the cut locus of $p$.
	
	Let $h= |\nabla \tilde f-\nabla e^u|$, $\Cal{F}_h(x,y) = \sup_{\gamma}\int_{\gamma}h $, where $\sup$
	is taken over all minimal geodesics $\gamma$ from $x$ to $y$. Let $\Psi=\Psi(\epsilon, L^{-1}|n,R,\rho)$.
	For $x_1\neq x_2 \in B_{(1-\Psi)\rho}(q)$, by Cheeger-Colding's segment inequality (\cite{Ch}, \cite{CC1}),
	$$\begin{aligned}
	&\int_{B_{\frac\Psi2}(x_1)\times B_{\frac\Psi2}(x_2)}\Cal{F}_h \\
	\leq &\; c(n,\rho)\left(\volume(B_{\frac{\Psi}{2}}(x_1))
	+ \volume(B_{\frac{\Psi}{2}}(x_2))\right)\int_{B_{\rho}(q)}|\nabla \tilde f-\nabla e^u|\\
	\leq &\; \Psi(\epsilon, L^{-1}|n,R,\rho).
	\end{aligned}$$
	Then there exists $x'_1 \in B_{\frac{\Psi}{2}}(x_1), x'_2 \in B_{\frac{\Psi}{2}}(x_2)$, such that
	$\int_{\gamma_{x'_1,x'_2}}h \leq \Psi(\epsilon, L^{-1}|n,R,\rho)$, i.e.,
	$$\left|\left(\tilde f(x'_1)- e^{u(x'_1)}\right) - \left(\tilde f(x'_2)- e^{u(x'_2)}\right)\right|\leq \Psi(\epsilon, L^{-1}|n,R,\rho). $$
	By Dirichlet-Poincar\'e inequality (\cite{Ch}), $$-\kern-1em\int_{B_{\rho}(q)}|\tilde f- e^u|\leq c(n, R)-
	\kern-1em\int_{B_{\rho}(q)}h \leq \Psi(\epsilon,L^{-1}|n,R,\rho).$$
	Consequently we obtain (5.3.2.1).
	
	Fixed $\alpha > 0$ small, by \cite{CC1} we can choose a cut-off function $\phi$ satisfying
	$$\begin{cases}
	\phi(x)= 1, & x\in B_{(1-\alpha)\rho}(q),\\
	\phi(x)= 0, & x \in  M\setminus B_{(1-\frac{\alpha}{2})\rho}(q),\end{cases} \quad |\nabla \phi|, |\Delta \phi|\leq c(n,\rho,\alpha).$$
	By (5.3.2.1), (5.3.2.2) and Bochner's formula, we derive
	$$\begin{aligned}
	&\Psi(\epsilon, L^{-1}|n,R,\rho,\alpha)\\
	\geq &\; \frac{1}{2}-\kern-1em\int_{B_{\rho}(q)}\Delta \phi (|\nabla \tilde f|^2 - \tilde f^2)\\
	= &\; \frac{1}{2}-\kern-1em\int_{B_{\rho}(q)}\phi \Delta(|\nabla \tilde f|^2 - \tilde f^2)\\
	= & -\kern-1em\int_{B_{\rho}(q)}\phi(|\op{Hess}\tilde f|^2 + \op{Ric}(\nabla \tilde f, \nabla \tilde f) + \langle\nabla\Delta \tilde f, \nabla \tilde f\rangle\\
	&  - \tilde f\Delta \tilde f - |\nabla \tilde f|^2)\\
	\geq & -\kern-1em\int_{B_{\rho}(q)}\phi(|\op{Hess}\tilde f - e^u|^2 +n e^{2u} - (n-1)e^{2u}\\
	& + n e^{2u}-n e^{2u} - e^{2u}) - \Psi(\epsilon, L^{-1}|n,R,\rho,\alpha)\\
	\geq & -\kern-1em\int_{B_{(1-\alpha)\rho}(q)}|\op{Hess}\tilde f - e^u|^2 - \Psi(\epsilon, L^{-1}|n,R,\rho,\alpha).
	\end{aligned}$$
\end{proof}

\end{document}